\input amstex\documentstyle{amsppt}  
\pagewidth{12.5cm}\pageheight{19cm}\magnification\magstep1
\topmatter
\title Conjugacy classes in reductive groups and two-sided cells \endtitle
\author G. Lusztig\endauthor
\address{Department of Mathematics, M.I.T., Cambridge, MA 02139}\endaddress
\thanks{Supported by NSF grant DMS-1566618.}\endthanks
\endtopmatter   
\document
\define\bV{\bar V}
\define\uR{\un R}
\define\uI{\un I}
\define\uV{\un V}
\define\uOm{\un\Om}

\define\uL{\un L}

\define\uW{\un W}

\define\mpb{\medpagebreak}

\define\bR{\bar R}

\define\frl{\forall}

\define\si{\sim}
\define\wt{\widetilde}
\define\sqc{\sqcup}

\define\qua{\quad}

\define\lb{\linebreak}

\define\op{\oplus}
   
\define\part{\partial}
\define\emp{\emptyset}

\define\ra{\rangle}
\define\n{\notin}
\define\iy{\infty}
\define\m{\mapsto}
\define\do{\dots}
\define\la{\langle}
\define\bsl{\backslash}

\define\lra{\leftrightarrow}

\define\sub{\subset}    

\define\T{\times}
\define\ti{\tilde}
\define\nl{\newline}
\redefine\i{^{-1}}

\define\un{\underline}
\define\ov{\overline}
\define\ot{\otimes}

\define\a{\alpha}
\redefine\b{\beta}
\redefine\c{\chi}
\define\g{\gamma}
\redefine\d{\delta}

\define\io{\iota}
\redefine\o{\omega}

\define\ph{\phi}
\define\ps{\psi}
\define\r{\rho}
\define\s{\sigma}

\define\th{\theta}

\redefine\l{\lambda}

\define\x{\xi}

\define\vt{\vartheta}

\redefine\G{\Gamma}
\redefine\D{\Delta}
\define\Om{\Omega}

\define\Ph{\Phi}

\define\boc{\bold c}

\define\kk{\bold k}

\redefine\ss{\bold s}

\define\BB{\bold B}
\define\CC{\bold C}

\define\NN{\bold N}

\define\QQ{\bold Q}

\define\TT{\bold T}

\define\ZZ{\bold Z}

\define\ca{\Cal A}
\define\cb{\Cal B}

\define\ch{\Cal H}
\define\ci{\Cal I}
\define\cj{\Cal J}

\define\cq{\Cal Q}
\define\car{\Cal R}
\define\cs{\Cal S}
\define\ct{\Cal T}
\define\cu{\Cal U}

\define\cw{\Cal W}

\define\fe{\frak e}

\define\fo{\frak o}

\define\fA{\frak A}

\define\fE{\frak E}

\define\fH{\frak H}

\define\fT{\frak T}
\define\fU{\frak U}

\define\te{\ti e}

\define\tm{\ti m}

\define\tw{\ti w}
\define\tz{\ti z}

\define\tT{\ti T}

\define\sha{\sharp}

\define\br{\bar r}

\define\bul{\bullet}

\define\BOU{B}
\define\MS{MS}
\define\SQINT{L1}
\define\ORA{L2}
\define\CELLIV{L3}
\define\RS{L4}
\define\HEC{L5}
\define\CDGVI{L6}
\define\CDGVII{L7}
\define\YOK{Y}

\head Introduction \endhead
\subhead 0.1\endsubhead 
Let $G'$ be a connected reductive group over $\CC$. In \cite{\CELLIV} it was shown that the set of
unipotent conjugacy classes of $G'$ is in natural bijection with the set of two-sided cells coming
from an affine Hecke algebra associated to $G$, a reductive group of type dual to that of $G'$. 
In this paper we extend this to a bijection between the set of all conjugacy classes of $G'$ and the
set of two-sided cells associated to a certain algebra $H$ containing the affine Hecke algebra, see Theorem 2.6.
The algebra $H$ is an affine analogue of an algebra considered in \cite{\CDGVI, \S31} which, on the one hand, is
a modified form of an algebra considered by Mars and Springer \cite{\MS} in their approach to the theory
of character sheaves and, on the other hand, as shown in \cite{\CDGVII, \S34}, is closely related to the
algebra studied by Yokonuma in \cite{\YOK}. Most of the paper is concerned with showing that the notion of
two-sided cell can be defined for $H$.

\subhead 0.2\endsubhead 
{\it Notation.}
Let $\car$ be a commutative ring with $1$ and let $A$ be an associative $\car$ algebra
with a given $\car$-basis $\{b_i;i\in I\}$ where $I\ne\emp$. For $K\sub I$ we write
$[K]=\sum_{i\in K}\car b_i\sub A$. For $i,i'\in I$ we write $i'\si_{left}i$ if for any 
$K\sub I$ such that $[K]$ is a left ideal of $A$ we have $i\in K$ if and only if $i'\in K$. 
This is an equivalence relation on $I$; the equivalence classes are called the left cells 
of $I$. Replacing ``left'' by ``right'' we obtain an equivalence relation $\si_{right}$ on 
$I$; the equivalence classes are called the right cells of $I$. Replacing ``left'' by 
``two-sided'' we obtain an equivalence relation $\si$ on $I$; the equivalence classes are 
called the two-sided cells of $I$.

By an affine Weyl group we mean a finite product of irreducible affine Weyl groups. In particular, the group with
one element is an affine Weyl group.

\head 1. Two-sided cells\endhead
\subhead 1.1\endsubhead
Let $\kk$, $\kk'$ be algebraically closed fields. Let $G$ be a connected reductive group 
over $\kk$ with a fixed Borel subgroup $\BB$ and a fixed maximal torus $\TT\sub\BB$. 
Let $\r=\dim\TT$, $\nu=\dim G/\BB$.

Let $G'$ be a connected reductive group over $\kk'$ of type dual to that of $G$ with a fixed 
Borel subgroup $\BB'$ and a fixed maximal torus $\TT'\sub\BB'$. 

Let $\uL$ (resp. $\uL'$) be the lattice of one parameter subgroups $\kk^*@>>>\TT$
(resp. $\kk'{}^*@>>>\TT'$). We assume that $\uL$ is also the group of characters
$\TT'@>>>\kk'{}^*$ and that $\uL'$ is also the group of characters $\TT@>>>\kk^*$.
We assume that the obvious nondegenerate pairing $\la,\ra:\uL\T\uL'@>>>\ZZ$ defined in 
terms of $\TT$ is the same as that defined in terms of $\TT'$.
Let $\uR\sub\uL$ (resp. $\uR'\sub\uL'$) be the set of coroots of $G$ (resp. $G'$) with 
respect to $\TT$ (resp. $\TT'$); we assume that $\uR$ (resp. $\uR'$) is also the set of
roots of $G'$ (resp. of $G$) with respect to $\TT'$ (resp. $\TT$). We also assume that the
canonical bijection $\uR\lra\uR'$ defined in terms of $G$ is the same as that defined in
terms of $G'$; we denote it by $h\lra h'$. Let $\{h_i;i\in\uI\}\sub\uR$ be the set of 
simple roots of $G'$ determined by $\BB'$; we assume that $\{h'_i;i\in\uI\}\sub\uR'$ is the 
set of simple roots of $G$ determined by $\BB$. 

Let $\uV=\QQ\ot\uL,\uV'=\QQ\ot\uL'$; now $\la,\ra$ extends to a nondegenerate bilinear
pairing $\la,\ra:\uV\T\uV'@>>>\QQ$. Let $V$ (resp. $V'$) be a $\QQ$-vector space containing 
$\uV$ (resp. $\uV'$) as a subspace and let $\la,\ra:V\T V'@>>>\QQ$ be a nondegenerate 
bilinear pairing extending $\la,\ra:\uV\T\uV'@>>>\QQ$.

We shall consider two cases (we will refer to them as case A and case B).

(A) We have $V=\uV$, $V'=\uV'$. We set $L=\uL,L'=\uL'$.

(B) We have $V=\uV\op\QQ c$, $V'=\uV'\op\QQ c'$ where $c,c'$ are vectors such that
$\la c,c'\ra=1$, $\la c,\uV'\ra=0$, $\la\uV,c'\ra=0$. We set $L=\uL\op\ZZ c$, 
$L'=\uL'\op\ZZ c'$. 

In case (A) we set $I=\uI$.

In case (B) we write $\uR=\sqc_{\fe\in\fE}\uR_\fe$, $\uR'=\sqc_{\fe\in\fE}\uR'_\fe$
where $(\uR_\fe,\uR'_\fe)$ are irreducible root systems and $\fE$ is an indexing set;
we have a corresponding partition $\uI=\sqc_{\fe\in\fE}\uI_\fe$ where 
$\uI_\fe$ indexes the simple roots/coroots in $\uR_\fe,\uR'_\fe$.
For each $\fe\in\fE$ let $\b_\fe\in\uR_\fe$ be such that
$-\b'_\fe$ is the highest root of $\uR'_\fe$. We set $I=\uI\sqc\fE$.
There is a unique function $\d:\uR@>>>\ZZ_{>0}$ such that $\d_{\b_\fe}=1$ for any $\fe\in\fE$
and for any $h,\ti h$ in $\uR$ we have $\la h,\ti h'\ra\d_{\ti h}=\la\ti h,h'\ra\d_h$. 
For any $h\in\uR$ we have $\d_h\in\{1,2,3\}$.
For $\fe\in\fE$ we set $h_\fe=\b_\fe+c\in L$, $h'_\fe=\b'_\fe\in L'$.

In both cases, for $i\in I$ we define a reflection $s_i:V@>>>V$ by 
$s_i(y)=y-\la y,h'_i\ra h_i$ and a reflection $s_i:V'@>>>V'$ by $s_i(x)=x-\la h_i,x\ra h'_i$
(this is the contragredient of $s_i:V@>>>V$). Let $W$ (resp. $\uW$) be the subgroup of 
$GL(V)$ generated by $\{s_i;i\in I\}$ (resp. by $\{s_i;i\in\uI\}$); taking contragredients we 
identify $W$ (resp. $\uW$) with the subgroup of $GL(V')$ generated by $\{s_i;i\in I\}$ 
(resp. by $\{s_i;i\in\uI\}$). It is well known that $W$ (resp $\uW$) is a Coxeter group on the set 
of generators $\{s_i;i\in I\}$ (resp. $\{s_i;i\in\uI\}$). In any case $\uW$ is a (finite) Weyl group.
In case A we have $W=\uW$; in case B, $W$ is an affine Weyl group (said to be the affine
Weyl group associated to $G$).

Let $w\m|w|$ be the length function $W@>>>\NN$; we set $S:=\{s_i;i\in I\}=\{w\in W;|w|=1\}$.
For $w\in W$ we have $w(L)\sub L$, $w(L')\sub L'$. 
Hence the $W$-action on $V'$ induces a $W$-action on $\bV':=V'/L'$. 
For $y\in L,\l\in\bV'$ we define $(y,\l)\in\QQ/\ZZ$ by $(y,\l)=\la y,x\ra$ where
$x\in V'$ is a representative of $\l$.

In case B we have $V'=\sqc_{\te\in\QQ}V'_{\te}$ 
where $V'_{\te}=\{x\in V';\la c,x\ra=\te\}$. The image of $V'_{\te}$ under 
$V'@>>>\bV'=V'/L'$ depends only on the image $e$ of $\te$ under the obvious map 
$\QQ@>>>\QQ/\ZZ$; we denote it by $\bV'_e$. We have $\bV'=\sqc_{e\in\QQ/\ZZ}\bV'_e$.
In case A we sometimes write $\bV'_0$ instead of $\bV'$. We show:

(a) {\it Any $W$-orbit on $\bV'$ is finite.}
\nl
In case A this is trivial since $W$ is finite. Assume now that we are in case B. Let 
$\l\in\bV'$ and let $x\in V'$ be a representative of $\l$. We can find an integer $q\ge1$ 
such that $qx\in L'$.
Let $\cq$ be the subgroup of $V'$ generated by $\{\d_hh';h\in\uR\}$; this is the free abelian
group with basis $\{\d_{h_i}h'_i;i\in\uI\}$. For $z\in\cq$ let $\th_z:V'@>>>V'$ be the
linear map $x\m x+\la c,x\ra z$; its contragredient is the linear map $y\m y-\la y,z\ra c$, 
$V@>>>V$. We have $\th_z\in W$ (compare \cite{\RS, 1.5}) and $z\m\th_z$ identifies $\cq$ 
with a normal subgroup $\cq'$ of $W$ such that $W=\cq'\uW$ (semidirect product). 

If $z\in\cq$ then $\th_{qz}(x)-x=\la c,qx\ra z\in\cq\sub L'$. Thus the stabilizer of $\l$ in
$W$ contains $\cq'{}^q$. Thus the $W$-orbit of $\l$ contains at most $\sha(W/(\cq'{}^q))$ 
elements. Since $\sha(W/\cq')<\iy$ and $\sha(\cq'/\cq'{}^q)<\iy$, we have 
$\sha(W/(\cq'{}^q))<\iy$. This proves (a).

\subhead 1.2\endsubhead
In case A we set $R=\uR$. In case B we set 
$$R=\{\c\in V;\c=w(h_i)\text{ for some $i\in I$ and some }w\in W\};$$
we have $R=\{h+\d_hmc;h\in\uR,m\in\ZZ\}$. In this case for $\c\in R$ we set $\c'=h'\in\uR'$ 
where $\c=h+\d_hmc$ with $h\in\uR,m\in\ZZ$. The map $R@>>>\uR'$, $\c\m\c'$ extends the
map $\uR@>>>\uR'$, $h\m h'$ considered earlier and it is compatible with the notation
$h'_\fe$ considered earlier.

In both cases, for $\c\in R$, the reflection $\s_\c:V@>>>V$, $y\m y-\la y,\c'\ra\c$ and its
contragredient $\s_\c:V'@>>>V'$, $x\m x-\la\c,x\ra\c'$ are defined and belong to $W$. Let 
$V^+=\sum_{i\in I}\QQ_{\ge0}h_i\sub V$. We set $R^+=R\cap V^+$, $R^-=-R^+$. 

Let $i_1,i_2,\do,i_r$ be a sequence in $I$ such that $w:=s_{i_1}s_{i_2}\do s_{i_r}$ 
satisfies $|w|=r$. A standard argument shows that

(a) {\it $\{\c\in R^+;w(\c)\in R^-\}$ consists of exactly $r$ elements, namely\lb
$s_{i_r}s_{i_{r-1}}\do s_{i_{j+1}}(h_{i_j})$ for $j=1,\do,r$.}

\mpb

For $\l\in\bV'$ we set $R_\l=\{\c\in R;(\c,\l)=0\in\QQ/\ZZ\}$, $R_\l^+=R_\l\cap R^+$, 
$R_\l^-=R_\l\cap R^-$; let $W'_\l=\{w\in W;w(\l)=\l\}$, a subgroup of $W$. Let $W_\l$ be the 
subgroup of $W$ generated by $\{s_\c;\c\in R_\l\}$. We have $W_\l\sub W'_\l$. If 
$w'\in W'_\l$ and $\c\in R_\l$, then $w'(\c)\in R_\l$; we have $w's_\c w'{}\i=s_{w'(\c)}$ 
hence $w'W_\l w'{}\i=W_\l$, so that $W_\l$ is a normal subgroup of $W'_\l$. 

In case A we set $\uR_\l=R_\l$. In case B we set 
$\uR_\l=\{h\in\uR;h+m\d_hc\in R_\l\text{ for some }m\in\ZZ\}$;
in the case where $\l\in\bV'_0$ we have $\uR_\l=\{h\in\uR;(h,\l)=0\in\QQ/\ZZ\}$.
In both cases we set $\uR'_\l=\{h';h\in\uR_\l\}\sub\uR'$. We show that for 
$h,\ti h\in\uR_\l$ we have ${}^1h:=s_{\ti h}(h)\in\uR_\l$. 

Assume first that we are in case B. We have $(h+m\d_hc,\l)=0$, $(\ti h+\tm\d_{\ti h}c,\l)=0$
for some $m,\tm$ in $\ZZ$. We have ${}^1h=h-\la h,\ti h'\ra\ti h\in\uR$,
$\d{{}^1h}=\d_h$. Let ${}^1m=m-\la\ti h,h'\ra\tm\in\ZZ$. We have
$$\align&({}^1h+{}^1m\d_{{}^1h}c,\l)=(h-\la h,\ti h'\ra\ti h+{}^1m\d_hc,\l)\\&=
-m\d_h(c,\l)-\la h,\ti h'\ra(\ti h,\l)+{}^1m\d_h(c,\l)\\&=
-m\d_h(c,\l)+\la h,\ti h'\ra\tm\d_{\ti h}(c,\l)+{}^1m\d_h(c,\l)\\&=
-m\d_h(c,\l)+\la \ti h,h'\ra\d_h(c,\l)+{}^1m\d_h(c,\l)=0.\endalign$$
Thus $({}^1h+{}^1m\d_{{}^1h}c,\l)=0$ so that ${}^1h\in\uR_\l$. 
The proof in case A is the same (we formally set $m=\tm={}^1m=0$ in the computation above.)
We see that $(\uR_\l,\uR_\l)$ is a root system. 

\mpb

We show:

(b) {\it For $\l\in\bV'$, $W_\l$ is a Coxeter group with length function 
$w\m|w|_\l=\sha(\c\in R_\l^+;w(\c)\in R_\l^-)$. Moreover, in case A, $W_\l$ is a Weyl group;
in case B, $W_\l$ is an affine Weyl group. If $s\in W_\l$, $|s|_\l=1$, then $s=s_\c$ where 
$\c$ is the unique element of $R_\l^+$ such that $s(\c)\in R_\l^-$.}
\nl
In case A, $W_\l$ is the Weyl group of $(\uR_\l,\uR'_\l)$ and (b) follows. We now assume 
that we are in case B. For any $\c\in R_\l$ let $\ch_\c=\{x\in V'_1;\la\c,x\ra=0\}$, a 
hyperplane in the affine space $V'_1$. For $w\in W_\l,\c\in R_\l$ we have 
$w(\ch_\c)=\ch_{w\i(\c)}$. For $x,x'$ in $V'_1-\cup_{\c\in R_\l}\ch_\c$ we say that 
$x\si x'$ if $\la\c,x\ra\la\c,x'\ra>0$ for any $\c\in R_\l$; this is an equivalence relation
on $V'_1-\cup_{\c\in R_\l}\ch_\c$; the equivalence classes are called $\l$-alcoves. There is
a unique $\l$-alcove $C_\l$ that contains $\{x\in V'_1;\la h_i,x\ra>0\frl i\in I\}$. Let 
$P_\l$ be the set of all $\c\in R_\l^+$  such that $\ch_\c$ is a wall of $C_\l$ (that is, 
the closure of $C_\l$ intersected with $\ch_\c$ is not contained in any affine subspace of 
$\ch_\c$ other than $\ch_\c$). We can apply the results in \cite{\BOU, Ch.V,\S3}, especially
Theorem 1; we see that $W_\l$ is a Coxeter group on the generators $\{s_\c;\c\in P_\l\}$ 
and with the length function as in (b). To prove the remaining statements of (b) we denote 
by $\uR_\l^u (u\in\cu)$ the irreducible components of $\uR_\l$. Note that if $u\ne u'$ in 
$\cu$, we have $\la h,\ti h'\ra=0$ for $h\in\uR_\l^u,\ti h\in\uR_\l^{u'}$. For $u\in\cu$ 
let $E_u$ be the affine subspace $c'+\sum_{h\in\uR_\l^u}\QQ h'$ of $V'_1$; let 
$R_\l^u=\{h+m\d_hc\in R_\l;h\in\uR_\l^u,m\in\ZZ\}$ and let $W_\l^u$ be the subgroup of 
$W_\l$ generated by $\{s_\c;\c\in R_\l^u\}$. Now for $\c\in R_\l^u$, the affine space $E_u$ 
is stable under $s_\c:V'@>>>V'$; hence $E_u$ is stable under $W_\l^u$. We have 
$W_\l=\prod_{u\in\cu}W_\l^u$. It is enough to show that for any $u\in\cu$, $W_\l^u$ is an 
affine Weyl group. For any $\c\in R_\l^u$ let 
$\ch_{\c,u}=\{x\in E_u;\la\c,x\ra=0\}=\ch_\c\cap E_u$, a hyperplane in $E_u$. For $x,x'$ in 
$E_u-\cup_{\c\in R_\l^u}\ch_{\c,u}$ we say that $x\si x'$ if $\la\c,x\ra\la\c,x'\ra>0$ for 
any $\c\in R_{\l,u}$; this is an equivalence relation on 
$E_u-\cup_{\c\in R_\l^u}\ch_{\c,u}$; the equivalence classes are called the alcoves of 
$E_u$. Using again the results in \cite{\BOU, Ch.V,\S3}, we see that it is enough to show 
that some (or equivalently any) alcove in $E_u$ is bounded. Let $h^j,j=1\do,t$ be the simple
coroots in $\uR_\l^u$ with respect to the set of positive coroots $\uR_\l^u\cap V^+$. Let 
$h^0$ be a coroot in $\uR_\l^u$ of the form $\sum_{j\in[1,t]}n_jh^j$ where $n_j\in\ZZ_{>0}$ 
for all $j$. For each 
$j\in[1,t]$ we can find $m_j\in\ZZ_{ge0}$ such that $\c^j:=h^j+m_j\d_jc\in R_\l^u$ (here
$\d_j=\d_{h^j}$). We can find $m_0\in\ZZ_{>0}$ such that $\c^0:=-h^0+m_0\d_0c\in R_\l^u$ 
(here $\d_0=\d_{h^0}$). Let 
$$\align&X=\{x\in E_u;\la\c^j,x\ra>0\frl j\in[0,t]\}\\&=\{x\in E_u;\la h^j,x\ra>-m_j\d_j
\frl j\in[1,t],\sum_{j\in[1,t]}n_j\la h^j,x\ra<m_0\d_0\}.\endalign$$
Note that $X$ contains at least one alcove and is bounded. (If $x\in E_u$, then 
$-m_j\d_j<\la h^j,x\ra<m_0\d_0/n_j$ for all $j\in[1,t]$.) It follows that some alcove in 
$E_u$ is bounded. This completes the proof of (b).

\mpb

We show:

(c) {\it Let $\l\in\bV'$. If $i\in I$ is such that $s_i\in W_\l$ then $h_i\in R_\l$.} 
\nl
We have $|s_i|_\l=\sha(\c\in R_\l^+;s_i(\c)\in R_\l^-)\le\sha(\c\in R^+;s_i(\c)\in R-)=1$. 
Thus $|s_i|_\l\le1$. Since $s_i$ is an element of $W_\l$ other than $1$ we have
$|s_i|_\l\ge1$ hence $|s_i|_\l=1$. By (b) we have $s_i=s_\c$ for a unique $\c\in R_\l^+$,
hence $\c=h_i$. This proves (c).

\mpb

Let $\l\in\bV'$. Let $i_1,i_2,\do,i_r$ be a sequence in $I$ such that 
$w:=s_{i_1}s_{i_2}\do s_{i_r}\in W'_\l$, $|w|=r$. Define 
$\ci=\{j\in[1,r];s_{i_r}\do s_{i_{j+1}}s_{i_j}s_{i_{j+1}}\do s_{i_r}\in W_\l$,
$X=\{\c\in R_\l^+;w(\c)\in R_\l^-\}$. We show:

(d) $\sha(\ci)=\sha(X)$.
\nl
Let $j\in\ci$. Then $s_{i_j}\in W_{s_{i_{j+1}}\do s_{i_r}(\l)}$ hence by (c) we have
$h_{i_j}\in R_{s_{i_{j+1}}\do s_{i_r}(\l)}$ so that 
$s_{i_r}\do s_{i_{j+1}}(h_{i_j})\in R_\l$. By (a), $s_{i_r}\do s_{i_{j+1}}(h_{i_j})$ 
($j\in\ci$) are distinct in $\{\c\in R^+;w(\c)\in R^-\}$. Thus,
 $j\m s_{i_r}\do s_{i_{j+1}}(h_{i_j})$ is an injective map $\ph:\ci@>>>X$. 

Now let $\c\in X$. From (a) we see that for some $j\in[1,r]$ we have
$\c=s_{i_r}\do s_{i_{j+1}}(h_{i_j})$. Since $\c\in R_\l$ we have
$$h_{i_j}\in  R_{s_{i_{j+1}}\do s_{i_r}(\l)}$$ so that
$$s_{i_j}\in W_{s_{i_{j+1}}\do s_{i_r}(\l)}$$ and therefore $j\in\ci$. We see that $\ph$ is
surjective hence a bijection. This proves (d).

\mpb

We show:

(e) {\it Any right coset $wW_\l\sub W$ ($w\in W$) contains a unique element $w_1$ such that
$w_1(R_\l^+)\sub R^+$. We have $|w_1|<|w_1z|$ for any $z\in W_\l-\{1\}$. We write 
$w_1=\min(wW_\l)$.}
\nl
In case A this is proved in \cite{\ORA, 1.9}. The following proof applies in both cases A,B.
We can find $w_1\in wW_\l$ such that $|w_1|\le|w_1z|$ for any $z\in W_\l$. We set $r=|w_1|$.
For any $\c\in R_\l^+$ we have $|w_1s_\c|\ne|w_1|$ hence $|w_1|<|w_1s_\c|$. 
Let $i_1,i_2,\do,i_r$ be a sequence in $I$ such that $s_\c w_1\i=s_{i_1}s_{i_2}\do s_{i_r}$.
$|s_\c w_1\i|=r$. Since $|s_\c(s_cw_1\i)|<|s_\c w_1\i|$, we see from \cite{\HEC, 2.2} that, 
for some $f\in[1,r]$ we have $s_\c=s_{i_1}s_{i_2}\do s_{i_f}\do s_{i_2}s_{i_1}$, hence
$\c=s_{i_1}s_{i_2}\do s_{i_{f-1}}(h_{i_f})$. Applying (a) with $w$ replaced by 
$w_1s_\c=s_{i_r}s_{i_{r-1}}\do s_{i_1}$, we see that
$w_1s_\c(\c)\in R^-$ that is $w_1(\c)\in R^+$. We have thus shown that $w_1(R_\l^+)\sub R^+$.
It follows that $w_1(R_\l^-)\sub R^-$. Now let $u\in W_\l-\{1\}$. We can find $\c\in R^+_\l$
such that $u(\c)\in R_\l^-$. We then have $w_1u(\c)\in R^-$. By the first part of the proof,
$w_1u$ is not of minimal length in $wW_\l$. Thus $|w_1|<|w_1u|$. This proves (e).

\mpb

Let $\fo$ be a $W$-orbit in $\bV'$. For $\l,\l'$ in $\fo$ let 
$$[\l',\l]=\{z\in W;\l'=z(\l),z=\min(zW_\l)\}=\{z\in W;\l'=z(\l),z(R^+_\l)=R^+_{\l'}\}.$$
Clearly, $[\l,\l']=[\l',\l]\i$ and for $\l''\in\fo$ we have 
$[\l'',\l'][\l',\l]\sub[\l'',\l]$. Thus, 

(f) {\it the group structure on $W$ makes $\G:=\sqc_{(\l',\l)\in\fo^2}[\l',\l]$ into a 
groupoid.}
\nl
In particular, for $\l\in\fo$, 
$$[\l,\l]=\{z\in W'_\l;z=\min(zW_\l)\}=\{z\in W'_\l;z(R^+_\l)=R^+_\l\}$$ 
is a subgroup of $W'_\l$. We show:

(g) {\it For $\l\in\bV'$, the group $W'_\l$ is the semidirect product of $[\l,\l]$ and 
$W_\l$ with $W_\l$ normal in $W'_\l$.}
\nl
By (e) we have $W'_\l=[\l,\l]W_\l$ and $[\l,\l]\cap W_\l=\{1\}$. It remains to recall that
$W_\l$ is normal in $W'_\l$.

\subhead 1.3\endsubhead
In this subsection we assume that we are in case B. We fix $e\in\QQ/\ZZ$, $\l\in\bV'_e$ and 
a representative $x$ for $\l$ in $V'$ such that $x\in V'_{\te}$ where $0\le\te<1$. We write 
$\te=p'/q'$ where $p',q'$ are integers, $q'\ge1$, $0\le p'<q'$ and $p',q'$ have no common 
prime divisor. 

\mpb

Assume first that $e=0$; then $\te=0$ so that $\la c,x\ra=0$. It follows that
$$R_\l=\{h+m\d_hc;h\in\uR_\l,m\in\ZZ\}.$$
We write $\uR_\l=\sqc_{u\in\cu}\uR_\l^u$ as in the proof of 1.2(b).
Let $h^{u,1},\do,h^{u,t_u}$ be the set of simple coroots of $\uR_\l^u$ with respect to the
set of positive coroots $\uR_\l^u\cap V^+$ and let $h^{u,0}$ be the coroot in $\uR_\l^u$ 
such that $(h^{u,0})'$ is minus the highest root of $(\uR_\l^u)'$. For $j\in[0,t_u]$, let 
$\d_{u,j}=\d_{h^{u,j}}$. 
Note that $h^{u,j}\in R^+_\l$ for $j\in[1,t_u]$ and $h^{u,0}+\d_{u,0}c\in R^+_\l$.
Let $\c\in R_\l^+$. We show that 
$$\c\in\sum_{u\in\cu,j\in[1,t_u]}\NN h^{u,j}+\sum_{u\in\cu}\NN(h^{u,0}+\d_{u,0}c).\tag a$$
We have $\c=h+m\d_hc$ where either $h\in\uR_\l^u,m\in\ZZ_{>0},u\in\cu$ or
$h\in\uR_\l^u\cap V^+,m=0,u\in\cu$. If $h\in\uR_\l^u\cap V^+,m=0,u\in\cu$, then clearly
$\c\in\sum_{j\in[1,t_u]}\NN h^{u,j}$. If $h\in\uR_\l^u,m\in\ZZ_{>0},u\in\cu$, then 
$$\align&\c=h+\d_hmc=(h-(\d_h/\d_{u,0})h^{u,0})+(\d_h/\d_{u,0})(m-1)(-h^{u,0})\\&
+(\d_h/\d_{u,0})m(h^{u,0}+\d_{u,0}c).\endalign$$
We now use $\d_h/\d_{u,0}\in\NN$, $h-(\d_h/\d_{u,0})h^{u,0}\in\sum_{j\in[1,t_u]}\NN h^{u,j}$
(a standard property of root systems) and (a) follows.
Using (a) and the proof of 1.2(b) we see that $W_\l$ is a Coxeter group on the generators 
$$s_{h^{u,j}}\text{ with }u\in\cu,j\in[1,t_u],\qua s_{h^{u,0}+\d_{u,0}c}\text{ with }
u\in\cu\}).\tag b$$

\mpb

Next we assume that $e\ne0$ (hence $0<p'<q'$) and $\uI=\{1\}$, $I=\uI\sqc\{i_0\}$, so that 
$\uR=\{h_1,-h_1\}$, $h_{i_0}=-h_1+c$, $h'_{i_0}=-h'_1$. We have
$x=(p/q)h'_1/2+(p'/q')c'$, where $p,q$ are integers, $q\ge1$ and $p,q$ have no common prime 
divisor.  We have
$$\align&R^+_\l=\{h_1+mc;m\in\ZZ_{\ge0},(p/q)+(p'/q')m\in\ZZ\}\\&
\sqc\{-h_1+mc;m\in\ZZ_{>0},-(p/q)+(p'/q')m\in\ZZ\}.\endalign$$
If $q'\n q\ZZ$ (so that $p\ne0$), then the equation $\pm(p/q)+(p'/q')m\in\ZZ$ has no integer 
solution $m$; hence in this case we have $R_\l=\emp$. We now assume that $q'\in q\ZZ$. We
have $\{m\in\ZZ;(p/q)+(p'/q')m\in\ZZ\}=m_1+q'\ZZ$,
$\{m\in\ZZ;-(p/q)+(p'/q')m\in\ZZ\}=m_2+q'\ZZ$, where $m_1,m_2$ are well defined integers
in $[1,q'-1]$ such that $p'(m_1+m_2)\in q'\ZZ$; thus we have $m_2=q'-m_1$ and
$$R^+_\l=\{h_1+mc;m=m_1+q'u,u\in\ZZ_{\ge0}\}\sqc\{-h_1+mc;m=q'-m_1+q'u,u\in\ZZ_{\ge0}\}.$$
In particular, we have $h_1+m_1c\in R^+_\l,-h_1+(q'-m_1)c\in R^+_\l$.
For $u\in\ZZ_{\ge0}$ we have

$h_1+(m_1+q'u)c=(u+1)(h_1+m_1c)+u(-h_1+(q'-m_1)c)$,

$-h_1+(q'-m_1+q'u)c=u(h_1+m_1c)+(u+1)(-h_1+(q'-m_1)c)$.
\nl
Thus $R_\l^+\sub\NN(h_1+m_1c)+\NN(-h_1+(q'-m_1)c)$. Using this and the proof of 1.2(b) we 
see that $W_\l$ is a Coxeter group (an infinite dihedral group) on the generators 
$s_{h_1+m_1c},s_{-h_1+(q'-m_1)c}$ (if $q'\in q\ZZ$) and $W_\l=\{1\}$ if $q'\n q\ZZ$.

\subhead 1.4\endsubhead
Until the end of 1.10 we fix a $W$-orbit $\fo$ in $\bV'$. 

Let $H_\fo$ be the $\ca$-algebra with $1$ defined  by
the generators $T_w (w\in W)$ and $1_\l (\l\in\fo)$ and relations:
$$1_\l1_\l=1_\l\text{ for }\l\in\bV',\qua 1_\l1_{\l'}=0\text{ for }\l\ne\l'\text{ in }\bV',$$
$$T_wT_{w'}=T_{ww'} \text{ for $ w,w'\in W$ such that }|ww'|=|w|+|w'|,$$ 
$$T_w1_\l=1_{w(\l)}T_w \text{ for }w\in W,\l\in\bV',$$
$$T_s^2=v^2T_1+(v^2-1)\sum_{\l\in\fo;s\in W_\l}T_s1_\l\text{ for }s\in S,$$
$$T_1=\sum_{\l\in\fo}1_\l.$$
Note that $T_1=\sum_{\l\in\fo}1_\l$ is the unit element of $H_\fo$. It will be convenient to
write $\tT_w=v^{-|w|}T_w$ for any $w\in W$. As in \cite{\CDGVI, 31.2}, we see that 

(a) {\it $\{\tT_w1_\l;w\in W,\l\in\fo\}$ is an $\ca$-basis of $H_\fo$.}
\nl
From the definitions we have:

(b) {\it if $s\in S,y\in W,\l\in\fo$ and $s\n W_{y(\l)}$, then
$\tT_\s\tT_y1_\l=\tT_{\s y}1_\l$.}
\nl
We show:

(c) {\it For $\l\in\fo$, $z\in[\l,\l],w\in W'_\l$ we have}
$\tT_z\tT_w1_\l=\tT_{zw}1_\l$, $1_\l\tT_{w\i}\tT_{z\i}=1_l\tT_{w\i z\i}$.
\nl
The proof is an extension of the proof of \cite{\CDGVII, 34.7(b)}. The use of the algebra 
antiautomorphism $\tT_x1_\l\m1_\l\tT_{x\i}$ shows that each of the two equalities in (c) 
implies the other equality. Let $i_1,i_2,\do,i_r$ be a sequence in $I$ such that 
$z\i=s_{i_1}s_{i_2}\do s_{i_r}$,
$|z|=r$. Applying 1.2(d) with $w$ replaced by $z\i$ and noting that $X$ in 1.2(d) is empty
in this case, we see that $s_{i_r}\do s_{i_{j+1}}s_{i_j}s_{i_{j+1}}\do s_{i_r}\n W_\l$ for
any $j\in[1,r]$. From the relations of $H_\fo$ we have
$$\tT_{w\i s_{i_1}s_{i_2}\do s_{i_{j-1}}}\tT_{s_{i_j}}1_{s_{i_{j+1}}\do s_{i_r}\l}=
\tT_{w\i s_{i_1}s_{i_2}\do s_{i_{j-1}}s_{i_j}}1_{s_{i_{j+1}}\do s_{i_r}\l}$$
for all $j\in[1,r]$. Using these identities we see that
$$\tT_{w\i}\tT_{z\i}1_\l=\tT_{w\i}\tT_{s_{i_1}}\tT_{s_{i_2}}\do \tT_{s_{i_r}}1_\l=
\tT_{w\i s_{i_1}s_{i_2}\do s_{i_r}}1_\l,$$
and (c) follows.

\mpb

The following result is a generalization of (c):

(d) {\it For $(\l',\l)in\fo^2,z\in[\l',\l],w\in W'_\l$ we have}
$\tT_z\tT_w1_\l=\tT_{zw}1_\l$, $1_\l\tT_{w\i}\tT_{z\i}=1_l\tT_{w\i z\i}$.
\nl
The use of the algebra antiautomorphism $\tT_x1_\l\m1_\l\tT_{x\i}$ shows that each of the 
two equalities in (d) implies the other equality. 
We can find $r\ge0$ and $i_1,i_2,\do,i_r$ in $I$ such that, setting 
$\l_0=\l,\l_1=s_{i_1}(\l),\l_2=s_{i_2}s_{i_1}(\l),\l_r=s_{i_r}\do s_{i_2}s_{i_1}(\l)$, we 
have $\l_0\ne\l_1\ne\l_2\ne\do\ne\l_r=\l'$.
We have $s_{i_1}\in[\l_1,\l_0],s_{i_2}\in[\l_2,\l_1],\do,s_{i_r}\in[\l_r,\l_{r-1}]$,
hence $s_{i_r}\do s_{i_2}s_{i_1}\in[\l_r,\l_0]=[\l',\l]$. We define $\tz\in W$ by 
$z=s_{i_r}\do s_{i_2}s_{i_1}\tz$. Then $\tz\in[\l,\l]$.
For $j\in[1,r]$ we have $s_{i_j}\n W_{s_{i_{j-1}}\do s_{i_1}\tz(\l)}$ 
(since $\l_i\ne\l_{i-1}$) hence, using (b) twice, we have
$$\tT_{s_{i_j}}\tT_{s_{i_{j-1}}\do s_{i_1}\tz}1_\l=\tT_{s_{i_j}s_{i_{j-1}}\do s_{i_1}1\tz}
1_\l,$$
$$\tT_{s_{i_j}}\tT_{s_{i_{j-1}}\do s_{i_1}\tz w}1_\l=\tT_{s_{i_j}s_{i_{j-1}}\do s_{i_1}
\tz w}1_\l.$$
Applying this repeatedly, we deduce 
$$\tT_{s_{i_r}\do s_{i_2} s_{i_1}\tz}1_\l=
\tT_{s_{i_r}}\do\tT_{s_{i_2}}\tT_{s_{i_1}}\tT_{\tz}1_\l,$$
$$\tT_{s_{i_r}\do s_{i_2}s_{i_1}\tz w}1_\l=\tT_{s_{i_r}}\do\tT_{s_{i_2}}\tT_{s_{i_1}}
\tT_{\tz w}1_\l.$$
By (c) we have $1_\l\tT_{w\i}\tT_{\tz\i}=1_\l\tT_{w\i\tz\i}$
and $\tT_{\tz w}1_\l=\tT_{\tz}\tT_w1_\l$. We deduce
$$\align&\tT_{zw}1_\l=\tT_{s_{i_r}\do s_{i_2}s_{i_1}\tz w}1_\l
=\tT_{s_{i_r}}\do\tT_{s_{i_2}}\tT_{s_{i_1}}\tT_{\tz w}1_\l\\&
=\tT_{s_{i_r}}\do\tT_{s_{i_2}}\tT_{s_{i_1}}\tT_{\tz}\tT_w1_\l
=\tT_{s_{i_r}}\do\tT_{s_{i_2}}\tT_{s_{i_1}}\tT_{\tz}1_\l\tT_w\\&
=\tT_{s_{i_r}\do s_{i_2}s_{i_1}\tz}1_\l\tT_w=\tT_z\tT_w1_\l.\endalign$$
This proves the first equality in (d); the second equality follows also.

\mpb

We show:

(e) {\it For $\l,\l',\l''$ in $\fo$, $z\in[\l',\l],z'\in[\l,\l'']$ we have}
$\tT_z\tT_{z'}1_{\l''}=\tT_{zz'}1_{\l''}$.
\nl
As in the proof of (d) we have $z=s_{i_r}\do s_{i_2}s_{i_1}\tz$,
$\tT_z1_\l=\tT_{s_{i_r}}\do\tT_{s_{i_2}}\tT_{s_{i_1}}\tT_{\tz}1_\l$,
where $i_1,\do i_r\in I$, $\l=\l_0\ne\l_1\ne\l_2\ne\do\ne\l_r=\l'$ are in $\bV'$,
$\l_0=\l,\l_1=s_{i_1}(\l),\l_2=s_{i_2}s_{i_1}(\l),\l_r=s_{i_r}\do s_{i_2}s_{i_1}(\l)$, 
$\tz\in[\l,\l]$.
Hence it is enough to show that
$$\tT_{s_{i_r}}\do\tT_{s_{i_2}}\tT_{s_{i_1}}\tT_{\tz}\tT_{z'}1_{\l''}=
\tT_{s_{i_r}\do s_{i_2}s_{i_1}\tz z'}1_{\l''}.$$
It is enough to show that for any $j\in[0,r]$ we have
$$\tT_{s_{i_j}}\do\tT_{s_{i_2}}\tT_{s_{i_1}}\tT_{\tz}\tT_{z'}1_{\l''}=
\tT_{s_{i_j}\do s_{i_2}s_{i_1}\tz z'}1_{\l''}.$$
We argue by induction on $j$. For $j=0$ we must show that
$\tT_{\tz}\tT_{z'}1_{\l''}=\tT_{\tz z'}1_{\l''}$; this follows from (d). Now assume that 
$j\ge1$. Using the induction hypothesis we see that it is enough to show that
$$\tT_{s_{i_j}}\tT_{w_j}1_{\l''}=\tT_{s_{i_j}w_j}1_{\l''}$$
where $w_j=s_{i_{j-1}}\do s_{i_2}s_{i_1}\tz z'$. Using (b) we see that it is enough to show 
that $s_{i_j}\n W_{s_{i_{j-1}}\do s_{i_2}s_{i_1}\tz z'(\l'')}$, or that $s_{i_j}\n W_{\l_j}$.
We have $s_{i_j}(\l_j)=\l_{j-1}\ne\l_j$ so that $s_{i_j}\n W_{\l_j}$ as desired.
This proves (e).

\mpb

We show:

(f) {\it Let $\l\in\fo$, $w\in W$, $\s\in W_\l$. Assume that $|\s|_\l=1$. We have
$\tT_w\tT_\s1_\l=\tT_{w\s}1_\l+\d(v-v\i)\tT_w1_\l$ where $\d\in\{0,1\}$. Moreover, if 
$w\in W_\l$ then $\d=0$ if $|w\s|_\l>|w|_\l$ and $\d=1$ if $|w\s|_\l<|w|_\l$.}
\nl
In case A this follows from \cite{\CDGVII, 34.7(a)}. The following proof applies in both
cases A,B. We have $\s=s_{i_1}s_{i_2}\do s_{i_r}$ for some sequence $i_1,i_2,\do,i_r$ in $I$
such that $r=|\s|$. From 1.2(d) we see that there is a unique $j\in[1,r]$ such that
$$s_{i_r}\do s_{i_{j+1}}s_{i_j}s_{i_{j+1}}\do s_{i_r}\in W_\l;$$ thus for $j'\in[1,r]-\{j\}$
we have $$s_{i_r}\do s_{i_{j'+1}}s_{i_{j'}}s_{i_{j'+1}}\do s_{i_r}\n W_\l;$$ moreover from the
proof of 1.2(d) we see that $s_{i_r}\do s_{i_{j+1}}(h_{i_j})$ is the unique element 
$\c\in R_\l^+$ such that $\s(\c)\in R_\l^-$. It follows that the reflection in $W$ defined
by $s_{i_r}\do s_{i_{j+1}}(h_{i_j})$ coincides with the reflection defined by 
$\c$; thus we have $$\s=s_{i_r}\do s_{i_{j+1}}s_{i_j}s_{i_{j+1}}\do s_{i_r}.$$
It follows that $$s_{i_1}s_{i_2}\do s_{i_{j-1}}s_{i_{j+1}}\do s_{i_r}=1.$$ From the relations
of $H_\fo$ we have
$$\align&\tT_{ws_{i_1}s_{i_2}\do s_{i_{j-1}}}\tT_{s_{i_j}}1_{s_{i_{j+1}}\do s_{i_r}(\l)}\\&=
\tT_{ws_{i_1}s_{i_2}\do s_{i_j}}1_{s_{i_{j+1}}\do s_{i_r}(\l)}+
\d'(v-v\i)\tT_{ws_{i_1}s_{i_2}\do s_{i_{j-1}}}1_{s_{i_{j+1}}\do s_{i_r}(\l)}\endalign$$
where $\d'=0$ if $|ws_{i_1}s_{i_2}\do s_{i_j}|>|ws_{i_1}s_{i_2}\do s_{i_{j-1}}|$ and 
$\d'=1$ otherwise,
$$\tT_{ws_{i_1}s_{i_2}\do s_{i_{j'-1}}}\tT_{s_{i_{j'}}}1_{s_{i_{j'+1}}\do s_{i_r}(\l)}=
\tT_{ws_{i_1}s_{i_2}\do s_{i_{j'}}}1_{s_{i_{j'+1}}\do s_{i_r}(\l)}$$
if $j'\in[1,r]-\{j\}$,
$$\align&\tT_{ws_{i_1}s_{i_2}\do s_{i_{j-1}}s_{i_{j+1}}\do s_{i_{j'-1}}}
\tT_{s_{i_{j'}}}1_{s_{i_{j'+1}}\do s_{i_r}(\l)}\\&=
\tT_{ws_{i_1}s_{i_2}\do s_{i_{j-1}}s_{i_{j+1}}\do s_{i_{j'}}}
1_{s_{i_{j'+1}}\do s_{i_r}(\l)}\endalign$$
if $j'\in[j+1,r]$. From these identities we see that
$$\align&\tT_w\tT_\s1_\l=\tT_w\tT_{s_{i_1}}\tT_{s_{i_2}}\do\tT_{s_{i_r}}1_\l\\&=
\tT_{ws_{i_1}s_{i_2}\do s_{i_{j-1}}}\tT_{s_{i_j}}\tT_{s_{i_{j+1}}}\do\tT_{s_{i_r}}1_\l\\&=
\tT_{ws_{i_1}s_{i_2}\do s_{i_{j-1}}s_{i_j}}\tT_{s_{i_{j+1}}}\do\tT_{s_{i_r}}1_\l\\&+
\d'(v-v\i)\tT_{ws_{i_1}s_{i_2}\do s_{i_{j-1}}}\tT_{s_{i_{j+1}}}\do\tT_{s_{i_r}}1_\l\\&=
\tT_{ws_{i_1}s_{i_2}\do s_{i_j}s_{i_{j+1}}\do s_{i_r}}1_\l\\&+
\d'(v-v\i)\tT_{ws_{i_1}s_{i_2}\do s_{i_{j-1}}s_{i_{j+1}}\do s_{i_r}}1_\l\\&=
=\tT_{w\s}1_\l+\d'(v-v\i)\tT_w1_\l.\endalign$$
Assume now that $w\in W_\l$. We show that $\d=\d'$. The condition that $\d=0$ is equivalent 
to the condition that $w(\c)\in R_\l^+$. The condition that $\d'=0$ is equivalent to the 
condition that $ws_{i_1}s_{i_2}\do s_{i_{j-1}}(h_{i_j})\in R^+$. It remains to note that
$\c=s_{i_1}s_{i_2}\do s_{i_{j-1}}(h_{i_j})$. (This follows from 
$\c=s_{i_r}\do s_{i_{j+1}}(h_{i_j})$ since
$s_{i_r}\do s_{i_{j+1}}=s_{i_1}s_{i_2}\do s_{i_{j-1}}$, or equivalently
$s_{i_1}s_{i_2}\do s_{i_{j-1}}s_{i_{j+1}}\do s_{i_r}=1$.)

\mpb

We show:

(g)  {\it Let $\l\in\fo$, $w\in W_\l,w'\in W_\l$. Assume that $|ww'|_\l=|w|_\l+|w'|_\l$. We 
have $\tT_{ww'}1_\l=(\tT_w1_\l)(\tT_{w'}1_\l)$.}
\nl
We argue by induction on $|w'|_\l$. If $w'=1$ the result is obvious. We now assume that 
$|w'|_\l\ge1$. We can write $w'=w'_1\s$ where $w'_1,\s$ are in $W_\l$ and 
$|w'_1|_\l=|w'|_\l-1$, $|\s|_\l=1$. We have $|ww'_1|_\l=|w|_\l+|w'_1|_\l$. From (f) we have 
$\tT_{w'}1_\l=(\tT_{w'_1}1_\l)(\tT_\s1_\l)$. Using this and the induction hypothesis we have
$$(\tT_w1_\l)(\tT_{w'}1_\l)=(\tT_w1_\l)(\tT_{w'_1}1_\l)(\tT_\s1_\l)=
(\tT_{ww'_1}1_\l)(\tT_\s1_\l).$$
It remains to use that $(\tT_{ww'_1}1_\l)(\tT_\s1_\l)=\tT_{ww'_1\s}1_\l$ which follows
again from (f).

\mpb

The following result is a special case of (f).

(h) {\it Let $\l\in\fo,\s\in W_\l$ be such that $|\s|_\l=1$. We have
$(\tT_\s1_\l)^2=1_\l+(v-v\i)\tT_\s1_\l$.}

\subhead 1.5\endsubhead
For $\l\in\fo$ let $H_\l$ be the $\ca$-submodule of $H_\fo$ with 
basis $\{\tT_w1_\l;w\in W_\l\}$. From 1.4(g),(h), we see that $H_\l$ is a subalgebra of 
$H_\fo$ (with a different unit, namely $1_\l$).

Let $\Xi=\{(\l',z,\l);(\l',\l)\in\fo^2,z\in[\l',\l]\}$. We show:

(a) {\it Let $(\l',z,\l)\in\Xi$. The assignment $w\m zwz\i$ is a Coxeter
group isomorphism $W_\l@>\si>>W_{\l'}$. The assignment $\tT_w1\l\m\tT_{zwz\i}1_{\l'}$
($w\in W_\l$) is an $\ca$-algebra isomorphism $\io_z:H_\l@>\si>>H_{\l'}$.}
\nl
Recall that $z(R^+_\l)=R^+_{\l'}$. This shows that $w\m zwz\i$ defines an isomorphism
$W_\l@>>>W_{\l'}$. For any $w\in W_\l$, $\c\m z(\c)$ defines a bijection
$\{\c_1\in R_\l^+;w(\c_1)\in R_\l^-\}@>>>\{\c_2\in R_{\l'}^+;zwz\i(\c_2)\in R_{\l'}^-\}$.
It follows that $|zwz\i|_{\l'}=|w|_\l$ for any $w\in W_\l$. Thus, $w\m zwz\i$ is
an isomorphism of Coxeter groups. We also see that $\tT_w1\l\m\tT_{zwz\i}1_{\l'}$ defines an
isomorphism of $\ca$-modules $\io_z:H_\l@>\si>>H_{\l'}$. To show that $\io_z$ is compatible 
with multiplication it is enough to note that, by 1.4(d),(e), for $w\in W_\l$ we have 
$\tT_{zwz\i}1_{\l'}=\tT_z\tT_w1_\l\tT_z\i$. This proves (a).

\mpb

Let $\car$ be a commutative ring with $1$. Let $Alg_\car$ be the category of associative
$\car$-algebras with $1$. Let $\Ph:\G@>>>Alg_\car$ be a functor; here the groupoid $\G$ (in 
1.2(f)) is viewed as a category in which all morphisms are isomorphisms. Now $\Ph$ consists
of a collection $\{\Ph(\l);\l\in\fo\}$ of associative $\car$-algebras with $1$ together 
with a collection of algebra isomorphisms $\io_z:\Ph(\l)@>\si>>\Ph(\l')$ (one for each 
$(\l',z,\l)\in\Xi$) such that for any $(\l'_1,z_1,\l_1),(\l'_2,z_2,\l_2)$ in $\Xi$ with
$\l'_2=\l_1$, the compositions 
$\Ph(\l_2)@>\io_{z_2}>>\Ph(\l'_2)=\Ph(\l_1)@>\io_{z_1}>>\Ph(\l'_1)$,
$\Ph(\l_2)@>\io_{z_1z_2}>>\Ph(\l'_1)$ coincide. 

Given $\Ph$ we define $\Ph^\bul\in Alg_\car$ as follows. As an $\car$-module we have 
$\Ph^\bul=\op_{(\l',z,\l)\in\Xi}\Ph_{\l',z,\l}^\bul$ where $\Ph_{\l',z,\l}^\bul=\Ph(\l)$. 
For any 
$(\l'_1,z_1,\l_1),(\l'_2,z_2,\l_2)$ in $\Xi$ and any 
$\x_1\in\Ph_{\l'_1,z_1,\l_1}^\bul=\Ph(\l_1)$, $\x_2\in\Ph_{\l'_2,z_2,\l_2}^\bul=\Ph(\l_2)$, 
we define the product $\x_1\x_2\in\Ph^\bul$ as $0$ if $\l'_2\ne\l_1$ and as 
$\io_{z_2}\i(\x_1)\x_2\in\Ph_{\l'_1,z_1z_2,\l_2}^\bul=\Ph(\l_2)$ (product in $\Ph(\l_2)$) if 
$\l'_2=\l_1$. From the definitions we see that this product is associative with unit element
such that its $\Ph_{\l',z,\l}^\bul$-component is 
$0$ unless $\l=\l'$, $z=1$. in which case it is the unit element of $\Ph(\l)$.

\subhead 1.6\endsubhead
Let $\ch:\G@>>>Alg_\ca$ be the functor defined by $\ch(\l)=H_\l$ for $\l\in\fo$ and the 
isomorphisms $\io_z:H_\l@>\si>>H_{\l'}$ in 1.5(a) for $(\l',z,\l)\in\Xi$. Then the 
associative algebra $\ch^\bul=\op_{(\l',z,\l)\in\Xi}\ch_{\l',z,\l}^\bul$ where 
$\ch_{\l',z,\l}^\bul=H_\l$ is defined as in 1.5.  We consider the $\ca$-linear map 
$\th:\ch^\bul@>>>H_\fo$ such that for any $(\l',z,\l)\in\Xi$ and any 
$\x\in\ch_{\l',z,\l}^\bul=H_\l$ we have $\th(\x)=\tT_z\x\in H_\fo$. From 1.2(e) we see that 
$\th$ is an isomorphism of $\ca$-modules. We show that it respects the algebra structures. 
For $(\l'_1,z_1,\l_1),(\l'_2,z_2,\l_2)$ in $\Xi$ and for $w_1\in W_{\l_1}$, 
$w_2\in W_{\l_2}$, we must show that 
$\tT_{z_1}\tT_{w_1}1_{\l_1}\tT_{z_2}\tT_{w_2}1_{\l_2}$ is zero if $\l_1\ne\l'_2$, while if 
$\l_1=\l'_2$, it is $\tT_{z_1z_2}\tT_{z_2\i w_1z_2}\tT_{w_2}1_{\l_2}$. The case where 
$\l_1\ne\l'_2$ is immediate. The case where $\l_1=\l'_2$ follows from 1.4(d),(e). We see 
that

(a) {\it $\th:\ch^\bul@>>>H_\fo$ is an isomorphism of $\ca$-algebras.}

\mpb

Now, for $\l\in\fo$, let $\fH_\l$ be the Hecke algebra associated to the Coxeter group 
$W_\l$. Thus, $\fH_\l$ is the $\ca$-module with basis $\{\fT_w;w\in W_\l\}$ with the 
associative $\ca$-algebra structure defined by the rules

$\fT_w\fT_{w'}=\fT_{ww'}$ if $w,w'$ in $W_\l$ satisfy $|ww'|_\l=|w|_\l+|w'|_\l$,

$\fT_\s^2=1+(v-v\i)\fT_\s$ if $\s\in W_\l$ satisfies $|\s|_\l=1$.
\nl
Note that $\fT_1$ is the unit element of $\fH_\l$. We define an $\ca$-linear map 
$\vt_\l:\fH_\l@>>>H_\l$ by $\fT_w\m\tT_w1_\l$ for $w\in W_\l$. This is clearly an 
isomorphism of $\ca$-modules; moreover, from 1.4(g),(h), we see that this is an algebra 
isomorphism.

Let $\ch':\G@>>>Alg_\ca$ be the functor defined by $\ch'(\l)=\fH_\l$ for $\l\in\fo$ and 
the isomorphisms $\io_z:\fH_\l@>\si>>\fH_{\l'}$, $\fT_w\m\fT_{zwz\i}$ for $w\in W_\l$ (here
$(\l',z,\l)\in\Xi$). Then the associative algebra 
$\ch'{}^\bul=\op_{(\l',z,\l)\in\Xi}\ch'_{\l',z,\l}{}^\bul$ where 
$\ch'_{\l',z,\l}{}^\bul=\fH_\l$ is defined as in 1.5.  

We consider the $\ca$-linear map $\vt:\ch'{}^\bul@>>>\ch^\bul$ such that for any 
$(\l',z,\l)\in\Xi$ and any $\x\in\ch'_{\l',z,\l}{}^\bul=\fH_\l$ we have 
$\vt(\x)=\vt_\l(\x)\in\ch_{\l',z,\l}^\bul=H_\l$. Clearly, $\vt$ is an isomorphism of 
$\ca$-modules. Moreover, $\vt$ is an isomorphism of algebras. (We use that for any 
$(\l',z,\l)\in\Xi$, the compositions $\fH_\l@>\vt_\l>>H_\l@>\io_z>>H_{\l'}$,
$\fH_\l@>\io_z>>\fH_{\l'}@>\vt_{\l'}>>H_{\l'}$, coincide.)

For $s\in S$ we have $\tT_s\i=\tT_s+(v\i-v)\sum_{\l\in\fo;s\in W_\l}1_\l$ in $H_\fo$. It 
follows that $\tT_w$ is invertible in $H_\fo$ for any $w\in W$. As in \cite{\CDGVI, 31.3},
there is a unique ring homomorphism $\bar{}:H_\fo@>>>H_\fo$ such that 
$\ov{\tT_w}=\tT_{w\i}\i$ for all $w\in W$ and $\ov{v^m1_\l}=v^{-m}1_\l$ for all $\l\in\fo$. 
The square of $\bar{}:H_\fo@>>>H_\fo$ is $1$. 

Now let $\l\in\fo$ and let $\s\in W_\l$ be such that $|\s|_\l=1$. From 1.4(h) we have
$\tT_\s\i1_\l=\tT_\s1_\l+(v\i-v)1_\l$ so that $\ov{\tT_\s1_\l}\in H_\l$. Using 1.4(g),(h) and
induction on $|w|_\l$, we see that for $w\in W_\l$ we have $\ov{\tT_w1_\l}\in H_\l$. Thus, 
$\bar{}:H_\fo@>>>H_\fo$ restricts to a ring isomorphism $\bar{}:H_\l@>>>H_\l$ whose square
is $1$. If $(\l',z,\l)\in\Xi$ then $\io_z:H_\l@>>>H_{\l'}$ is compatible with the bar 
operators on $H_\l$, $H_{\l'}$. (It is enough to check that for any $\s\in W_\l$ such that 
$|\s|_\l=1$ we have $\io_z(\tT_\s\i1_\l)=\tT_{z\s z\i}\i1_{\l'}$; this is immediate.) It 
follows that the group homomorphism $\bar{}:\ch^\bul@>>>\ch^\bul$ whose restriction to any 
$\ch_{\l',z,\l}^\bul=H_\l$ is $\bar{}:H_\l@>>>H_\l$ is in fact a ring homomorphism with square 
$1$. We show that $\th:\ch^\bul@>>> H_\fo$ is compatible with $\bar{}:\ch^\bul@>>>\ch^\bul$ and 
$\bar{}:H_\fo@>>>H_\fo$. It is enough to show that for $(\l',z,\l)\in\Xi$, $w\in W_\l$, we 
have $\ov{\tT_z\tT_w1_\l}=\tT_z\tT_{w\i}\i1_\l$; it is also enough to show that 
$\ov{\tT_z1_\l}=\tT_z1_\l$ or that $\tT_{z\i}\i1_\l=\tT_z1_\l$, which follows from 1.4(e).

For $\l\in\fo$, $\bar{}:H_\l@>>>H_\l$ corresponds under $\vt_\l$ to a ring isomorphism
$\bar{}:\fH_\l@>>>\fH_\l$ such that $\ov{v^m\fT_w}=v^{-m}\fT_{w\i}\i$ for any $w\in W_\l$,
$m\in\ZZ$. Under the isomorphism $\vt:\ch'{}^\bul@>>>\ch^\bul$, $\bar{}:\ch^\bul@>>>\ch^\bul$ 
corresponds to a ring isomorphism $\bar{}:\ch'{}^\bul@>>>\ch'{}^\bul$ with square $1$; its 
restriction to $\ch'_{\l',z,\l}{}^\bul=\fH_\l$ coincides with $\bar{}:\fH_\l@>>>\fH_\l$.

Let $\l\in\fo$. By \cite{\HEC, 5.12}, for any $w\in W_\l$ there is a unique element 
$c_w^\l\in\fH_\l$ such that $\ov{c_w^\l}=c_w^\l$ and 
$c_w^\l-\fT_w\in\sum_{y\in W_\l;|y|_\l<|w|_\l}v\i\ZZ[v\i]\fT_y$.
Using $\vt_\l$ we deduce that for any $w\in W_\l$ there is a unique element 
$c_{w,\l}\in H_\l$ such that $\ov{c_{w,\l}}=c_{w,\l}$ and 
$c_{w,\l}-\tT_w1_\l\in\sum_{y\in W_\l;|y|_\l<|w|_\l}v\i\ZZ[v\i]\tT_y1_\l$.

If $(\l',z,\l)\in\Xi$ and $w\in W_\l$ then $\tT_zc_{w,\l}\tT_{z\i}$ satisfies the
definition $c_{zwz\i,\l'}\in H_{\l'}$ hence we have 
$\tT_zc_{w,\l}\tT_{z\i}=c_{zwz\i,\l'}$.

Now let $w\in W,\l\in\fo$. We can write uniquely $w=z\tw$ where $\tw\in W_\l$, 
$z\in[\l',\l]$ for some
$\l'\in\fo$. We set $c_{w,\l}=\tT_zc_{\tw,\l}$. We have $\ov{c_{w,\l}}=c_{w,\l}$ 
and $c_{w,\l}-\tT_w1_\l\in\sum_{y\in W_\l;|y|_\l<|\tw|_\l}v\i\ZZ[v\i]\tT_{zy}1_\l$. It 
follows that $\{c_{w,\l};w\in W,\l\in\fo\}$ is an $\ca$-basis of $H_\fo$.

\subhead 1.7\endsubhead
For any $\l\in\fo$ we state some properties of the multiplication of elements of form 
$c_{w,\l}\in H_\l$ with $w\in W_\l$. These properties follow from the known properties of 
the multiplication of elements of form $c_w^\l\in\fH_\l$ via the isomorphism $\vt_\l$. (We 
use that $W_\l$ is a Weyl group or an affine Weyl groups, see 1.2(b).)
For any $w,w'$ in $W_\l$ we have  
$$c_{w,\l}c_{w',\l}=\sum_{w''\in W_\l}r^{w'';\l}_{w,w'}c_{w'',\l}$$
where $r^{w''}_{w,w'}\in\ca$ are zero for all but finitely many $w''$. Moreover, for any 
$w''\in W_\l$ there is a well defined smallest integer $a(w'')\ge0$ such that
$v^{a(w'')}r^{w'';\l}_{w;w'}\in\ZZ[v]$ for all $w,w'$ in $W_\l$. We have 
$$r^{w'';\l}_{w,w'}=\br^{w'';\l}_{w;w'}v^{-a(w'')}\mod v^{-a(w'')+1}\ZZ[v]$$
where $\br^{w'';\l}_{w,w'}\in\NN$.

We now consider the free abelian group $H_\l^\iy$ with basis $\{t_w;w\in W_\l\}$. Then
$H_\l^\iy$ is an (associative) ring with multiplication given by
$$t_wt_{w'}=\sum_{w''\in W_\l}\br^{w'';\l}_{w,w'}t_{w''}.$$
It has a unit element of the form $1=\sum_{w\in D_\l}t_w$ where $D_\l$ is a well defined 
finite subset of $\{w\in W_\l;w^2=1\}$. The $\ca$-linear map 
$\ps_\l:H_\l@>>>\ca\ot H_\l^\iy$, given by
$$\ps_\l(c_{w,\l})=\sum_{w''\in W_\l,d\in D_\l;a(d)=a(w'')}r^{w'';\l}_{w,d}t_{w''}$$
is a homomorphism of $\ca$-algebras with $1$.

For $w,w'$ in $W_\l$ we have $w\si_{left}w'$ (resp. $w\si_{right}w'$ if and only if for 
some $u\in W_\l$, $t_{w'}$ appears with non-zero coefficient in $t_ut_w$ (resp. in $t_ut_w$).
For $w,w'$ in $W_\l$ we have $w\si w'$ if and only if for some $u,u'$ in $W_\l$ we have 
$t_{w'}$ appears with non-zero coefficient in $t_ut_wt_{u'}$ (or equivalently if for some 
$w''\in W_\l$ we have $w\si_{left}w''$, $w'\si_{right}w''$). Recall 
from 0.2 that the equivalence classes for $\si_{left}$ (resp. $\si_{right}$, $\si$) are 
called the left (resp. right, two-sided) cells of $W_\l$. For any left (resp. right, 
two-sided) cell $K$ of $W_\l$ the subgroup $\sum_{w\in K}\ZZ t_w$ is a left (resp. right, 
two-sided) ideal of $H_\l^\iy$. 

Note that the definition of left (resp. right, two-sided) cells of $W_\l$ depends only on
the Coxeter group structure of $W_\l$ hence these are defined for any Coxeter group $\cw$
which is a Weyl group or an affine Weyl group; we shall write $Cell(\cw)$ for the set of
two-sided cells of $\cw$.

\mpb

Let $\cj:\G@>>>Alg_\ZZ$ be the functor defined by $\cj(\l)=H^\iy_\l$ for $\l\in\fo$ and the 
isomorphisms $\io_z:H^\iy_\l@>\si>>H^\iy_{\l'}$ for $(\l',z,\l)\in\Xi$ where $\io_z$ maps 
any basis element $t_w (w\in W_\l)$ to the basis element $t_{zwz\i}$ of $H^\iy_{\l'}$. Then 
the (associative) ring $\cj^\bul$ is defined as in 1.5.

\subhead 1.8\endsubhead
For $(w_0,\l_0),(w_1,\l_1)$ in $W\T\fo$ we write in $H_\fo$:
$$c_{w_0,\l_0}c_{w_1,\l_1}=\sum_{(w_2,\l_2)\in W\T\fo}R^{w_2,\l_2}_{w_0,\l_0;w_1,\l_1}
c_{w_2,\l_2}$$
where $R^{w_2,\l_2}_{w_0,\l_0;w_1,\l_1}\in\ca$ are zero for all but finitely many 
$(w_2,\l_2)$. For $j=0,1,2$ we can write uniquely $w_j=z_j\tw_j$ where $z_j\in[\l'_j,\l_j]$,
$\tw_j\in W_{\l_j}$. We have 
$$c_{w_0,\l_0}c_{w_1,\l_1}=\tT_{z_0}c_{\tw_0,\l_0}\tT_{z_1}c_{\tw_1,\l_1}.$$
This is $0$ unless $\l_0=z_1(\l_1)$. Assume now that $\l_0=z_1(\l_1)$. We have
$$\align&c_{w_0,\l_0}c_{w_1,\l_1}=\tT_{z_0}\tT_{z_1}c_{z_1\i\tw_0z_1,\l_1}c_{\tw_1,\l_1}\\&=
\tT_{z_0z_1}\sum_{y\in W_{\l_1}}r_{z_1\i\tw_0z_1,\l_1;\tw_1,\l_1}^{y,\l_1}c_{y,\l_1}\\&=
\sum_{y\in W_{\l_1}}r_{z_1\i\tw_0z_1,\l_1;\tw_1,\l_1}^{y,\l_1}c_{z_0z_1y,\l_1}.\endalign$$
We see that $R_{w_0,\l_0;w_1,\l_1}^{w_2,\l_2}=0$ unless the conditions 
$\l_2=\l_1,\l_0=z_1(\l_1),z_2=z_0z_1$ are satisfied; if these conditions are satisfied, then
$$R_{w_0,\l_0;w_1,\l_1}^{w_2,\l_2}=r_{z_1\i\tw_0z_1,\tw_1}^{\tw_2;\l_1}.$$
We see that

(a) $v^{a_{w_2,\l_2}}R_{w_0,\l_0;w_1,\l_1}^{w_2,\l_2}\in\ZZ[v]$
\nl
where $a_{w_2,\l_2}=a_{\tw_2}$ is defined in terms of $W_{\l_2}$ and that $a_{w_2,\l_2}$ is 
the smallest integer such that (a) holds for any $(w_0,\l_0),(w_1,\l_1)$. We now define 
$\bR_{w_0,\l_0;w_1,\l_1}^{w_2,\l_2}\in\ZZ$ by
$$R_{w_0,\l_0;w_1,\l_1}^{w_2,\l_2}=\bR_{w_0,\l_0;w_1,\l_1}^{w_2,\l_2}v^{-a(w_2,\l_2)}
\mod v^{-a(w_2,\l_2)+1}\ZZ[v].$$
We consider the free abelian group $\fH^\iy_\fo$ with basis
$\{t_{w,\l};w\in W,\l\in\fo\}$. We define a ring structure on $\fH^\iy_\fo$ by 
$$t_{w_0,\l_0}t_{w_1,\l_1}=\sum_{(w_2,\l_2)\in W\T\fo}\bR^{w_2,\l_2}_{w_0,\l_0;w_1,\l_1}
t_{w_2,\l_2},$$
or equivalently (for $\tw_0\in W_{\l_0},\tw_1\in W_{\l_1}$, $z_0\in[\l'_0,\l_0]$,
$z_1\in[\l'_1,\l_1]$):
$$t_{z_0\tw_0,\l_0}t_{z_1\tw_1,\l_1}=\sum_{\tw_2\in W_{\l_1}}
\br^{\tw_2;\l_1}_{z_1\i\tw_0z_1,\tw_1}t_{z_0z_1\tw_2,\l_1},\tag b$$
if $\l_0=z_1(\l_1)$,
$$t_{z_0\tw_0,\l_0}t_{z_1\tw_1,\l_1}=0, \text{ if }\l_0\ne z_1(\l_1).\tag c$$
This ring is associative. Indeed, this ring has the same multiplication rule as the ring 
$\cj^\bul$ in 1.7, which is known to be associative from 1.5. Note that the ring $\fH^\iy_\fo$
has a unit element, namely $\sum_{\l\in\fo}\sum_{w\in D_\l}t_{w,\l}$.

The $\ca$-algebra homomorphisms $\ps_\l:H_\l@>>>\ca\ot H_\l^\iy$ combine to give a 
$\ca$-algebra homomorphism $\ps:\ch^\bul@>>>\ca\ot\cj^\bul$ or equivalently 
$\ps:H_\fo@>>>\ca\ot\fH^\iy_\fo$; it is given by
$$\ps(c_{w,\l})=\sum_{(w_2,\l_2)\in W\T\fo,\l_1\in\fo,d\in D_{\l_1};a(w_2,\l_2)=a(d,\l_1)}
R^{w_2,\l_2}_{w,\l;d,\l_1}t_{w_2,\l_2}$$
that is,
$$\ps(c_{w,\l})=\sum_{(\tw_2\in W_\l,d\in D_\l;a(\tw_2)=a(d)}r^{\tw_2;\l}_{\tw,d}
t_{z\tw_2,\l}$$
where $w=z\tw$ with $\tw\in W_\l$, $z\in[w(\l),\l]$.

\subhead 1.9\endsubhead
We now describe the left (resp. right) cells of $W\T\fo$
defined in terms of the basis $(t_{w,\l})$ of $H_\fo^\iy$. (See 0.2.) 
Let $(w,\l),(w',\l')$ in $W\T\fo$. We say that 
$(w,\l)\approx_{left}(w',\l')$ (resp. $(w,\l)\approx_{right}(w',\l')$) if for some 
$(u,\l_1)\in W\T\fo$, $t_{w',\l'}$ appears with non-zero coefficient in the product
$t_{u,\l_1}t_{w,\l}$ (resp. $t_{w,\l}t_{u,\l_1}$). 

We write $w=z\tw$, $w'=z'\tw'$ where $\tw\in W_\l$, $\tw'\in W_{\l'}$, $z\in[w(\l),\l]$,
$z'\in[w'(\l'),\l']$. From 1.8(b),(c) we see that the condition that
$(w,\l)\approx_{left}(w',\l')$ is that $\l=\l'$ and $\tw\si_{left}\tw'$ in $W_\l=W_{\l'}$;
the condition that $(w,\l)\approx_{right}(w',\l')$ is that $\l=(z\i z')(\l')$ and 
$\tw'\si_{right}(z'{}\i z)\tw(z\i z')$ in $W_{\l'}=W_{(z'{}\i z)(\l)}$.

Using the results in 1.7 we deduce that $\approx_{left}$ and $\approx_{right}$ are 
equivalence relations on $W\T\fo$. In particular, $\approx_{left}$ is transitive, hence if 
for $(w,\l)\in W\T\fo$ we set
$\wt{(w,\l)}=\{(w',\l')\in W\T\fo;(w,\l)\approx_{left}(w',\l')\}$, then the subgroup spanned 
by $\{t_{w',\l'};(w',\l')\in\wt{(w,\l)}\}$ is a left ideal of $H^\iy_\fo$, so that
$\wt{(w,\l)}$ is a union of left cells. 
If $(w,\l)\approx_{left}(w',\l')$, then clearly any left ideal of $H^\iy_\fo$ spanned by 
a subset of the canonical basis of $H^\iy_\fo$ and containing $t_{w,\l}$ must also
contain $t_{w',\l'}$. Since $\approx_{left}$ is symmetric we must also have
$(w',\l')\approx_{left}(w,\l)$ hence any left ideal of $H^\iy_\fo$ spanned by a subset of 
the canonical basis of $H^\iy_\fo$ and containing $t_{w',\l'}$ must also contain $t_{w,\l}$;
it follows that $(w,\l),(w',\l')$ are in the same left cell. We now see that any equivalence
class for $\approx_{left}$ is exactly one left cell. We see that the canonical basis
elements indexed by any left cell span a left ideal. The same argument shows that the 
equivalence classes for $\approx_{right}$ are exactly the right cells and that the canonical
basis elements indexed by any right cell span a right ideal.

We now see that we have a natural bijection
$$\{\text{left cells in }W\T\fo\}\lra\{(\l,\g);\l\in\fo,\g=\text{left cell in }W_\l\}$$
The bijection associates to each $(\l,\g)$ in the right hand side the left cell
$\{(z\tw,\l);\tw\in\g, z=\min(zW_\l)\}$ of $W\T\fo$. 
We also see that the right cells in $W\T\fo$ are the images of the left cells under the
involution $W\T\fo@>>>W\T\fo$, $(w,\l)\m(w\i,w(\l))$.

\subhead 1.10\endsubhead
We now describe the two-sided cells of $W\T\fo$
defined by the basis $(t_{w,\l})$ of $H_\fo^\iy$. (See 0.2). 
Let $\l,\l'$ in $\fo$; let $w,w'$ in $W$.

We say that $(w,\l)\approx(w',\l')$ if for some $(u_1,\l_1),(u_2,\l_2)$ in $W\T\fo$,
$t_{w',\l'}$ appears with non-zero coefficient in the product
$t_{u_1,\l_1}t_{w,\l}t_{u_2,\l_2}$ or equivalently, if for some $(w'',\l'')$ we have
$(w,\l)\approx_{left}(w'',\l'')$ and $(w',\l')\approx_{right}(w'',\l'')$.
(The equivalence uses the positivity of the structure constants of $H^\iy_\fo$.)
We write $w=z\tw$, $w'=z'\tw'$ where $\tw\in W_\l$, $\tw'\in W_{\l'}$, $z\in[w(\l),\l]$,
$z'\in[w'(\l'),\l']$. From 1.9 we see that the condition that $(w,\l)\approx(w',\l')$ is 
that for some $\l''\in\fo$, $\tw''\in W_{\l''}$, $z_1\in[\l',\l'']$ we have $\l=\l''$, 
$\tw\si_{left}\tw''$ in $W_\l=W_{\l''}$, $\tw''\si_{right}z_1\i\tw'z_1$ in 
$W_{\l''}=W_{z_1\i(\l')}$ or equivalently: 
for some $z_1\in[\l',\l]$ we have $\tw\si z_1\i\tw'z_1$ in $W_\l=W_{z_1\i(\l')}$. 
Using the results in 1.7 we deduce that $\approx$ is an equivalence relation on $W\T\fo$.

In particular, $\approx$ is transitive, hence if for $(w,\l)\in W\T\fo$ we set
$\hat{(w,\l)}=\{(w',\l')\in W\T\fo;(w,\l)\approx(w',\l')\}$, then the subgroup spanned 
by $\{t_{w',\l'};(w',\l')\in\hat{(w,\l)}\}$ is a two-sided ideal of $H^\iy_\fo$, so that
$\hat{(w,\l)}$ is a union of two-sided cells. 

Now assume that $(w,\l)\approx(w',\l')$; let $(w'',\l'')$ be such that
$(w,\l)\approx_{left}(w'',\l'')$ and $(w',\l')\approx_{right}(w'',\l'')$. By 1.9, we have 
$(w,\l)\si_{left}(w'',\l'')$ and $(w',\l')\si_{right}(w'',\l'')$. Hence if $\ci$ is a 
two-sided ideal of $H^\iy_\fo$ spanned by a subset of the canonical basis then we have 
$t_{w,\l}\in\ci\lra t_{w'',\l''}\in\ci\lra t_{w',\l'}\in\ci$. It follows that 
$(w,\l)\si(w',\l')$. We see that any equivalence class for $\approx$ is contained in a 
two-sided cell. As we have seen earlier, any equivalence class for $\approx$ is a union of 
two-sided cells, hence it is exactly one two-sided cell. Also, the canonical basis elements 
indexed by any equivalence class for $\approx$ (hence by any two-sided cell) span a 
two-sided ideal of $H^\iy_\fo$.

\mpb

Let $Cell(W\T\fo)$ be the set of two-sided cells of $W\T\fo$.
Note that for any $\l\in\fo$, the group $[\l,\l]$ acts by conjugation on $W_\l$; this
action induces an action of $[\l,\l]$ on $Cell(W_\l)$. If $\D\in Cell(W\T\fo)$, 
for any $\l\in\fo$ we set 
$\D_\l=\{\tw\in W_\l;(\tw,\l)\in\D\}$; this is a union of the two-sided cells of $W_\l$ in a 
fixed $[\l,\l]$-orbit. Moreover, if $\l,\l'$ are elements of $\fo$ then
$z\D_\l z\i=D_{\l'}$ for any $z\in[\l',\l]$. We now see that for each $\l\in\fo$ we have a 
bijection $\D\m\D_\l$ between $Cell(W\T\fo)$ and the set of orbits 
of the conjugation of $[\l,\l]$ on $Cell(W_\l)$. 

From the results in 1.8 we see that if $\D,\l$ are as above, then the function 
$(w,\l)\m a(w,\l)$ on $\D$ is constant with value equal to the value of the $a$-function
on any of the two-sided cells of $W_\l$ contained in $\D_\l$.

\subhead 1.11\endsubhead
Let $H=\op_\fo H_\fo$ (resp. $H^\iy=\op_\fo H^\iy_\fo$) where $\fo$ runs over all $W$-orbits
in $\bV'_0$. We view $H$ (resp. $H^\iy$) as an $\ca$-algebra (resp. ring) without $1$ in 
general) in which $H_\fo$ (resp. $H^\iy_\fo$) is a subalgebra (resp. subring) for any $\fo$ 
and $H_\fo H_{\fo'}=0$ (resp. $H^\iy_\fo H^\iy_{\fo'}=0$) for $\fo\ne\fo'$.
Note that $\{c_{w,\l};w\in W,\l\in\bV'_0\}$ is an $\ca$-basis of $H$ and
 that $\{t_{w,\l};w\in W,\l\in\bV'_0\}$ is a $\ZZ$-basis of $H^\iy$.
The left (resp. right, two-sided) cells of $W\T\bV'_0$ are defined in terms of the basis
$(t_{w,\l})$ of $H^\iy$ as in 0.2.
Using 1.9 we have a natural bijection
$$\{\text{left cells in }W\T\bV'_0\}\lra\{(\l,\g);\l\in\bV'_0,\g=\text{left cell in }W_\l\}$$
The bijection associates to each $(\l,\g)$ in the right hand side the left cell
$\{(z\tw,\l);\tw\in\g, z=\min(zW_\l)\}$ of $W\T\bV'_0$. 
The right cells in $W\T\bV'_0$ are the images of the left cells under the
involution $W\T\bV'_0@>>>W\T\bV'_0$, $(w,\l)\m(w\i,w(\l))$.
Let $Cell(W\T\bV'_0)$ be the set of two-sided cells of $W\T\bV'_0$. 
We have $Cell(W\T\bV'_0)=\sqc_\fo Cell(W\T\fo)$ 
where $\fo$ runs over all $W$-orbits in $\bV'_0$ and $Cell(W\T\fo)$ is described in 1.10. 

Let $(w'',\l'')\in W\T\bV'_0$. There is a well defined smallest integer $a(w,\l)\ge0$
such that for any $(w,\l),(w',\l')$ in $W\T\bV'_0$, the coefficient of $c_{w'',\l''}$ in 
the product $c_{w,\l}c_{w',\l'}$ belongs to $v^{-a(w'',\l'')}\ZZ[v]$.
From 1.10 we see that for any $\D\in Cell(W\T\bV'_0)$ the function $\D@>>>\NN$,
$(w'',\l'')\m a(w'',\l'')$ is constant on $\D$.

\head 2. Conjugacy classes\endhead
\subhead 2.1\endsubhead
Let $\cw$ be an affine Weyl group. Let $\ct$ be the set of translations in $\cw$, that is
the elements of $\cw$ whose $\cw$-conjugacy class is finite. Note that $\ct$ is a subgroup
of finite index in $\cw$. Let $\cs$ be the set of simple reflections of $\cw$. Let $A(\cw)$
be the group of automorphisms $\s:\cw@>>>\cw$ such that $\s(\cs)=\cs$; note that 
$\s\in A(\cw)$ is uniquely determined by its restriction to $\cs$.

We now assume that the affine Weyl group $\cw$ is irreducible. Let 
$\cs_0$ be the set of all $s\in\cs$ such 
that $\cs-\{s\}$ together with $\ct$ generates $\cw$. We have $\cs_0\ne\emp$. 
If $\s\in A(\cw)$ then $\s$ restricts to a permutation of $\cs_0$. Let $A_0(\cw)$
be the set of all $\s\in A(\cw)$ such that $\s:\cs_0@>>>\cs_0$ is either fixed point free 
or the identity. Note that $A_0(\cw)$ is a normal subgroup of $A(\cw)$; it acts simply
transitively on $\cs_0$. For any $s\in\cs_0$ let $A^s(\cw)$ be the set of all $\s\in A(\cw)$
such that $\s(s)=s$ (a subgroup of $A(\cw)$); we have $A(\cw)=A_0(\cw)A^s(\cw)$, 
$A_0(\cw)\cap A^s(\cw)=\{1\}$.

\subhead 2.2\endsubhead
The results in this subsection can be deduced from those in 2.1.
Let $\cw$ be an affine Weyl group with set of simple reflections
$\cs$. We can write $\cw=\prod_{u\in\cu}\cw^u$, $\cs=\sqc_{u\in\cu}\cs^u$ where
$\cw^u (u\in\cu)$ is a finite collection of irreducible affine Weyl groups and $\cs^u$ is
the set of simple reflections of $\cw^u$. We set $\cs_0=\sqc_{u\in\cu}\cs^u_0$ where 
$\cs^u_0\sub\cs^u$ is as in 2.1 (with $\cw,\cs$ replaced by $\cw^u,\cs^u$). We define 
$A(\cw)$ as in 2.1. Let $A_0(\cw)$ be the set of all $\s\in A(\cw)$ such that for any 
$u\in\cu$, we have $\s(\cw^u)=\cw^u$ and $\s|_{\cw^u}\in A_0(\cw^u)$. 
Let $\ss$ be a subset of $\cs_0$ such that for any $u\in\cu$, $\ss$ contains exactly one 
element of $\cs^u_0$. (Note that $A_0(\cw)$ acts simply transitively on the set of all such
$\ss$.) Let $A^\ss(\cw)$ be the set of all $\s\in A(\cw)$ such that 
$\s(\ss)=\ss$ (a subgroup of $A(\cw)$); we have $A(\cw)=A_0(\cw)A^\ss(\cw)$, 
$A_0(\cw)\cap A^\ss(\cw)=\{1\}$.

\subhead 2.3\endsubhead
In the setup of 2.2 let $\s\in A(\cw)$. Let $\ct$ be as in 2.1. The following property 
(which is checked case by case) has been stated in \cite{\SQINT}:

(a) {\it we have $\s\in A_0(\cw)$ if and only if there exists $w\in\cw$ such that 
$\s(t)=wtw\i$ for all $w\in W$.}
\nl
Note that the image of $w$ in $\cw/\ct$ is uniquely determined by $\s$. In particular, 
$\s\m w\ct$ is a well defined (injective) homomorphism $A_0(\cw)@>>>\cw/\ct$. We show:

(b) {\it If $\s\in A_0(\cw)$ and $c\in Cell(\cw)$, then $\s(c)=c$.}
\nl
This is proved in \cite{\CELLIV, 4.9(b)} using the property (a) of $\s$.

\subhead 2.4\endsubhead
In this subsection we assume that we are in case B. Let $\l\in\bV'_0$. Recall that 
$\uR_\l=\{h\in\uR;(h,\l)=0\in\QQ/\ZZ\}$, $R_\l=\{h+m\d_hc\in R;h\in\uR_\l,m\in\ZZ\}$.
Let $\uR_\l^u (u\in\cu)$ be the irreducible components of $\uR_\l$.
Note that $\uR_\l^+=\uR_\l\cap V^+$ is a set of positive (co)roots for $\uR_\l$ with set of 
simple (co)roots $\un\Pi_\l=\{h^{u,j};u\in\cu,j\in[1,t_u]\}\sub\uR^+_\l$ (notation of 1.3). 
Let ${}'\un\Pi_\l=\{h^{u,0};u\in\cu\}\sub R_\l$ (notation of 1.3),
${}'\Pi_\l=\{h^{u,0}+\d_{h^{u,0}}c;u\in\cu\}\sub R^+_\l$,
$\hat\Pi_\l=\un\Pi_\l\sqc{}'\un\Pi_\l\sub R_\l$,
$\Pi_\l=\un\Pi_\l\sqc{}'\Pi_\l\sub R^+_\l$.
Note that we have a unique bijection $\Pi_\l@>\si>>\hat\Pi_\l$, $\c\m\hat\c$, such that for 
$\c\in\Pi_\l$ we have $\hat\c\in\c+\ZZ c$. 
Recall from 1.3 that $W_\l$ is an affine Weyl group with Coxeter generators
$\cs_\l=\{s_\c;\c\in\Pi_\l\}$. The results of 2.2, 2.3 are applicable with 
$(\cw,\cs,\ct)=(W_\l,\cs_\l,\cq)$.
Let $\uW_\l$ be the subgroup of $\uW$ generated by $\{s_\c;\c\in\uR_\l\}$ (a subgroup of 
$W_\l$. Note that $\uW_\l$ is a Weyl group with Coxeter generators $\{s_h;h\in\un\Pi_\l\}$.
Let $\uW'_\l=\uW\cap W'_\l=\{w\in\uW;w(\l)=\l\}$. We have $\uW_\l\sub\uW'_\l$. We write 
$\Om_\l$ instead of
$$[\l,\l]=\{z\in W'_\l;z(R^+_\l)=R^+_\l\}=\{z\in W'_\l;z(\Pi_\l)=\Pi_\l\}.$$
Let
$$\uOm_\l=\{z\in\uW'_\l;z(\uR^+_\l)=\uR^+_\l\}=\{z\in\uW'_\l;z(\un\Pi_\l)=\un\Pi_\l\}.$$
We have $\uW'_\l=\uOm_\l\uW_\l$ and $\uOm_\l\cap\uW_\l=\{1\}$. We show:

(a) $\uOm_\l\sub\Om_\l$.
\nl
Let $z\in\uOm_\l$. Let $u\in\cu$. For some $u'\in\cu$, $h\m z(h)$ defines an isomorphism
$\uR_\l^u@>\si>>\uR_\l^{u'}$ which takes simple roots of $\uR_\l^u$ to simple roots of 
$\uR_\l^{u'}$; hence it takes $h^{u,0}$ to $h^{u',0}$ and we must have 
$\d_{h^{u,0}}=\d_{h^{u',0}}$. We have also $z(c)=c$, since $z\in\uW$. Thus we have 
$z(h^{u,0}+\d_{h^{u,0}}c)=h^{u',0}+\d_{h^{u',0}}c$. We see that $z({}'\Pi_\l)\sub{}'\Pi_\l$,
hence $z({}'\Pi_\l)={}'\Pi_\l$. We see that $z(\Pi_\l)=\Pi_\l$ so that $z\in\Om_\l$. This 
proves (a).

\mpb

Let $\fA_\l$ be the group of permutations $\r:\Pi_\l@>>>\Pi_\l$ with the following property:
there exists (a necessarily unique) $\s=\s_\r\in A(W_\l)$ such that $\s(s_\c)=s_{\r(\c)}$ 
for any $\c\in\Pi_\l$. Note that $\r\m\s_\r$ is an isomorphism $\fA_\l@>\si>>A(W_\l)$. Each 
$\r\in\fA_\l$ defines a permutation 
$\hat\r:\hat\Pi_\l@>>>\hat\Pi_\l$ by $\hat\r(\hat\c)=\widehat{\r(\c)}$ for any 
$\c\in\Pi_\l$. Now $\r\m\hat\r$ defines an isomorphism of $\fA_\l$ onto a subgroup
$\hat\fA_\l$ of the group of permutations of $\hat\Pi_\l$.
Let $\ti\Om_\l$ be the set of all $z\in\uW'_\l$ such that $z(\hat\Pi_\l)=\hat\Pi_\l$ 
and the permutation of $\hat\Pi_\l$ defined by $h\m z(h)$ belongs to $\hat\fA_\l$.

Note that any element $w\in W$ can be written uniquely in the form $w=t\un w$ where 
$\un w\in\uW$, $t\in\cq$; moreover $w\m\un w$ is a homomorphism $\ti k:W@>>>\uW$ with
kernel $\cq$. We show:

(b) {\it  we have $\ti k(\Om_\l)\sub\ti\Om_\l$ hence $\ti k$ defines a homomorphism 
$k:\Om_\l@>>>\ti\Om_\l$ with kernel contained in $\cq$.}
\nl
Let $w\in\Om_\l$. We write $w=t\un w$ with $\un w\in\uW$, $t\in\cq$. Since $\l\in\bV'_0$, we
have $t(\l)=\l$. Thus $\un w(\l)=t\i w(\l)=t\i(\l)=\l$, so that $\un w\in\uW'_\l$. Now let 
$\c\in\Pi_\l$ so that $\hat\c\in\hat\Pi_\l$ and $w(\c)\in\Pi_\l$. We show that 
$\un w(\hat\c)=\widehat{w(\c)}$. Now $\un w(\hat\c)\in\uR$ and $h=\widehat{w(\c)}$ is the 
unique element in $\uR$ such that $w(\c)-h\in\ZZ c$. Thus it is enough to show that
$w(\c)-\un w(\hat\c)\in\ZZ c$. We have $\c=\hat\c+mc$ with $m\in\ZZ$ hence it is enough to
show that $w(\hat\c+mc)-\un w(\hat\c)\in\ZZ c$. Since $\un w(c)=c$ and $w=t\un w$, it is 
enough to show that $t\un w(\hat\c)-\un w(\hat\c)\in\ZZ c$ or, setting 
$h'=\un w(\hat\c)\in\uR$, that $t(h')-h'\in\ZZ c$. The last identity holds for any 
$t\in\cq,h'\in\uR$, as we can see using the definitions. This proves (b).

\mpb

We shall need the following property:

(c) {\it Let $\r\in\fA_\l$ be such that $\s_\r\in A_0(W_\l)$; let $\hat\r$ be the 
corresponding element of $\hat\fA_\l$. Then there is a unique $w\in\uW_\l$ such that 
$w(h)=\hat\r(h)$ for any $h\in\hat\Pi_\l$.}
\nl
Let $w\in\uW$ be such that $\s_\r(t)=wtw\i$ for all $t\in\cq$. (We use 2.3(a).) A case by 
case check shows that $w$ satisfies the requirement of (c). The uniqueness of $w$ is 
obvious. 

\mpb

We define a homomorphism 
$\a:\Om_\l@>>>fA_\l$ by $z\m\r$ where $\r(\c)=z(\c)$ for any $\c\in\Pi_\l$.
We define a homomorphism $\un\a:\uOm_\l@>>>\fA_\l$ as the composition 
$\uOm_\l@>>>\Om_\l@>\a>>\fA_\l$ (the first map is the inclusion (a)). It follows that

(d) $\text{image}(\un\a)\sub\text{image}(\a)$.
\nl
We define a homomorphism $\ti\a:\ti\Om_\l@>>>\fA_\l$ as the composition
$\ti\Om_\l@>>>\hat\fA_\l@>>>\fA_\l$ where the first homomorphism is $z\m\r'$ where
$\r'(h)=z(h)$ for any $h\in\hat\Pi_\l$ and the second homomorphism is $\hat\r@>>>\r$
for $\r\in\fA_\l$. We show that $\a=\ti\a k$ with $k$ as in (b).
Let $w\in\Om_\l$. We write $w=t\un w$ where $\un w\in\uW$, $t\in\cq$.
We must show that for $\c\in\Pi_\l$ we have $\widehat{w(\c)}=\un w(\hat\c)$. This has been
verified in the course of proving (b). From $\a=\ti\a k$ it follows that

(e) $\text{image}(\a)\sub\text{image}(\ti\a)$.

\mpb

We define $\a':A_0(W_\l)@>>>\fA_\l$ as the composition $A_0(W_\l)@>>>A(W_\l)@>>>\fA_\l$
where the first map is the obvious imbedding and the second map is the inverse of the 
bijection $\r\m\s_\r$. We show:

(f) {\it We have $\text{image}(\ti\a)\sub\text{image}(\a')\text{image}(\un\a)$.}
\nl
Let $z\in\ti\Om_\l$. Let $\r=\ti\a(z)\in\fA_\l$ so that $\hat\r(h)=z(h)$ for any 
$h\in\hat\Pi_\l$. Let $\s=\s_\r\in A(W_\l)$. By 2.2 we can write $\s=\s'\s''$ where 
$\s'\in A_0(W_\l)$ and $\s''\in A(W_\l)$ maps $\{s_\c;\c\in\Pi_\l-\un\Pi_\l\}$ into itself. 
Let $\r',\r''$ in $\fA_\l$ be such that $\s'=\s_{\r'},\s''=\s_{\r''}$. Then $\r=\r'\r''$ 
and $\r''$ maps $\Pi_\l-\un\Pi_\l$ into itself and $\un\Pi_\l$ into itself. We have 
$\hat\r=\hat\r'\hat\r''$ and $\hat\r''$ maps $\un\Pi_\l$ into itself (recall that for 
$\c\in\un\Pi_\l$ we have $\hat\c=\c$). By (c), we can find $w\in\uW_\l$ such that 
$w(h)=\hat\r'(h)$ for any $h\in\hat\Pi_\l$. For any $h\in\hat\Pi_\l$ we have 
$\hat\r'\hat\r''(h)=z(h)$ hence $w(\hat\r''(h))=z(h)$, that is $w\i z(h)=\hat\r''(h)$. 
Since $\hat\r''(\un\Pi_\l)=\un\Pi_\l$, we deduce that $w\i z(\un\Pi_\l)=\un\Pi_\l$. Since 
$z\in\un W'_\l$, $w\in\un W_\l$, we have $w\i z\in\un W'_\l$. We see that $w\i z\in\uOm_\l$,
so that $z=wx$ where $x\in\uOm_\l$. We have $w\in\ti\Om_\l,x\in\ti\Om_\l$ hence
$\r=\ti\a(z)=\ti\a(w)\ti\a(x)$. We have $\ti\a(w)\in\text{image}(\a')$,
$\ti\a(x)\in\text{image}(\un\a)$. This proves (f).

\mpb

We show:

(g) {\it Let $c\in Cell(W_\l)$. The collection of two-sided cells $\{zcz\i;z\in\Om_\l\}$ 
coincides with the collection of two-sided cells $\{zcz\i;z\in\uOm_\l\}$.}
\nl
From (e),(f) we deduce that $\text{image}(\a)\sub\text{image}(\a')\text{image}(\un\a)$.
Hence, if $z\in\Om_\l$ then the automorphism $\o\m z\o z\i$ of $W_\l$ is a product $\s'\s''$
where $\s'\in A_0(W_\l)$ and $\s''$ is conjugation by an element $z'\in\uOm_\l$. Thus, we 
have $zcz\i=\s'(z'cz'{}\i)$. By 2.3(b) we have $\s'(z'cz'{}\i)=z'cz'{}\i$ hence
$zcz\i=z'cz'{}\i$. We see that the first collection in (g) is contained in the second
collection. The reverse containment follows from (a). This proves (g).

\subhead 2.5\endsubhead
In the remainder of this section we assume that $\kk'$ has characteristic zero.
Let $\TT'_f$ (resp. $\kk'_f{}^*$) be the set of elements of finite order in $\TT'$ (resp. 
$\kk'{}^*$). We have canonically $\TT'=\uL'\ot\kk'{}^*$, $\TT'_f=\uL'\ot\kk'{}^*_f$. 
Recall that in case B we have $\bV'_0=\uV'/\uL'$, see 1.3, and that in case A we write 
$\bV'_0=\bV'$.
We have canonically $\bV'_0=\uL'\ot(\QQ/ZZ)$.
Note that $\uW$ can be viewed as the Weyl group of $G$ and that of $G'$; it acts naturally
on $\TT$ and $\TT'$. We choose an isomorphism $\QQ/\ZZ@>\si>>\kk'_f{}^*$.
Via this isomorphism we have $\TT'_f=\uL'\ot\QQ/\ZZ=\bV'_0$.
This is compatible with the $\uW$-actions on $\TT'$ and $\bV'_0$.
Now let $\fo$ be a $W$-orbit on $\bV'$ which is contained in $\bV'_0$. 
Note that $\fo$ is in fact a $\uW$-orbit. In case A this is because $\uW=w$; in case B this 
is because the subgroup $\cq$ (see the proof of 1.1(a)) of $W$ acts trivially on $\bV'_0$.
 Under the identification
$\TT'_f=\bV'_0$, $\fo$ can be viewed as a $\uW$-orbit on $\TT_f$.
Let $\l\in\fo$. Let $Z(\l)$ be the centralizer of $\l\in\TT'$ in $G'$ and let $Z(\l)^0$ be
the identity component of $Z(\l)$. Since $\l$ is semisimple in $G'$, $Z(\l)^0$ is a reductive
group and $\TT'$ is a maximal torus of it.
Let $h\in\uR$. We have $h\in\uL$ so that $h$ can be viewed as a homomorphism $\TT'@>>>K^*$.
The condition that $h$ is a root of $Z(\l)^0$ with respect to $\TT'$ is that $h(\l)=1$. An
equivalent condition (with $\l$ is viewed as an element of $\bV'_0$) is that $(h,\l)=0$ or 
that $h\in\uR_\l$. We see that the set of roots of $Z(\l)^0$ is $\uR_\l$ and the 
corresponding set of coroots is $\uR'_\l\sub\uR'$ (as in 1.2). Using this and 1.3(b), we 
see that we have $W_\l=\cw(Z(\l)^0)$ where $\cw(Z(\l)^0)$ is the Weyl group (in case A) or 
the affine Weyl group (in case B) associated to the dual of the reductive group of 
$Z(\l)^0$, in the same way as $W$ is the Weyl group (in case A) or the affine Weyl group 
(in case B) associated to $G$, the dual of $G'$.
Let $\fU(Z(\l)^0)$ be the set of unipotent conjugacy classes of $Z(\l)^0$; let
$\fU_{sp}(Z(\l)^0)$ be the set of special unipotent conjugacy classes of 
$Z(\l)^0$ (a subset of $\fU(Z(\l)^0)$). By \cite{\CELLIV, 4.8(b)} in case B we have a 
canonical bijection 

(a) $\fU(Z(\l)^0)\lra Cell(\cw(Z(\l)^0))=Cell(W_\l)$;
\nl
this restricts to a bijection $\fU_{sp}(Z(\l)^0)\lra Cell(\uW_\l)$ (each two-sided cell of
$\uW_\l$ is contained in a two-sided cell of $W_\l$ and this gives an imbedding
$Cell(\uW_l)\sub Cell(W_\l)$. The last bijection can be viewed as a bijection

(b) $\fU_{sp}(Z(\l)^0)\lra Cell(W_\l)$ 
\nl
in case A.

Now let $\fU(Z(\l))$ be the set of unipotent elements of $Z(\l)^0$ up to $Z(\l)$-conjugacy;
let $\fU_{sp}(Z(\l))$ be the set of special unipotent elements of $Z(\l)^0$ up to 
$Z(\l)$-conjugacy. 
Let $N(\TT')$ be the normalizer of $\TT'$ in $G'$. If $w\in\uW'_\l$ and $n\in N(\TT')$ is a 
representative of $w$ viewed as an element of $N(\TT')/\TT'$ then $n\in Z(\l)$ and
$u\m nun\i$ defines  bijections $\fU(Z(\l)^0)@>>>\fU(Z(\l)^0)$,
$\fU_{sp}(Z(\l)^0)@>>>\fU_{sp}(Z(\l)^0)$ which depend only on $w$ not on $n$ (since
$\TT'\sub Z(\l)^0$). This gives an action of $\uW'_\l$ on $\fU(Z(\l)^0)$ leaving stable 
$\fU_{sp}(Z(\l)^0)$. (In this action the subgroup $\uW_\l$ acts trivially.)
It is easy to see that two unipotent $Z(\l)^0$-conjugacy classes of 
$Z(\l)^0$ are in the same $\uW'_\l$-orbit if and only if they are contained in the same
$Z(\l)$-conjugacy class. 
Let $\uOm_\l$ be as in 2.4 (in case B) and let $\uOm_\l=[\l,\l]$ in case A.
Since $\uW'_\l=\uOm_\l\uW_\l$ it follows 
that two unipotent $Z(\l)^0$-conjugacy classes of 
$Z(\l)^0$ are in the same $\uOm_\l$-orbit if and only if they are contained in the same
$Z(\l)$-conjugacy class. 
Thus we can identify $\fU(Z(\l))=\uOm_\l\bsl\fU(Z(\l)^0)$,
$\fU_{sp}(Z(\l))=\uOm_\l\bsl\fU_{sp}(Z(\l)^0)$.
From the definitions we see that the bijections (a),(b) are compatible with the natural
actions of $\uW'_\l$ on the two sides of (a),(b). Taking orbits of these actions we deduce
bijections

(c) $\fU(Z(\l))\lra \uOm_\l\bsl Cell(W_\l)$ (in case B);

(d) $\fU_{sp}(Z(\l))\lra [\l,\l]\bsl Cell(W_\l)$ (in case A).
\nl
Now in case B the natural action of $\uOm_\l$ on $Cell(W_\l)$ extends to an action of 
$[\l,\l$ on $Cell(W_\l)$ with the same orbit space (see 2.4(g)). Thus (c) can be viewed as
a bijection

(e) $\fU(Z(\l))\lra [\l,\l]\bsl Cell(W_\l)$ (in case B).
\nl
By 1.10 we have a canonical bijection $Cell(W\T\fo)\lra[\l,\l]\bsl Cell(W_\l)$ (in both
cases A,B). Combining with (d),(e) we obtain bijections

$Cell(W\T\fo)\lra\fU(Z(\l))$ in case B;

$Cell(W\T\fo)\lra\fU_{sp}(Z(\l))$ in case A.
\nl
Using this and the equality $Cell(W\T\bV'_0)=\sqc_\fo Cell(W\T\fo)$ in 1.11 we see that
$Cell(W\T\bV'_0)$ is in natural bijection with $\sqc_\l\fU(Z(\l))$ (in case B) and with 
$\sqc_\l\fU_{sp}(Z(\l))$ (in case A) where $\l$ runs though a set of representatives for 
the $\uW$-orbits in $\bV'_0=\TT'$. 
Let $Conj(G')$ be the set of conjugacy classes of elements in $G'$
with semisimple part of finite order. Let $Conj_{sp}(G')$ be the set of conjugacy classes 
of elements in $G'$ with semisimple part of finite order and such that the unipotent
part is special in the connected centralizer of the semisimple part. Note that we have
canonically $\sqc_\l\fU(Z(\l))=Conj(G')$, $\sqc_\l\fU_{sp}(Z(\l))=Conj_{sp}(G')$. 
For any $g\in G'$ let $\cb'_g$ be the variety of Borel
subgroups of $G'$ containing $g$. We have the following result.

\proclaim{Theorem 2.6} (a) We have canonical bijections $Cell(W\T\bV'_0)\lra Conj_{sp}(G')$ 
(in case A) and $Cell(W\T\bV'_0)\lra Conj(G')$ (in case B).

(b) If $\D\in Cell(W\T\bV'_0)$ and $\boc$ is the conjugacy class in $G'$ corresponding to
$\D$ under (a), then the value of the $a$-function (see 1.11) on $\D$ is equal to 
$\dim\cb'_g$ for $g\in\boc$.
\endproclaim
Now (a) is obtained by combining several statements above; (b) follows from 1.10, 1.11,
using \cite{\CELLIV, 4.8(c)}.

\enddocument